\documentclass[a4paper,11pt,reqno]{amsart}

\usepackage[headings]{fullpage}

\usepackage{dsfont,amssymb,graphicx}
\usepackage[all]{xy}

\usepackage{hyperref}

\usepackage[nobysame,alphabetic,initials,msc-links]{amsrefs}

\DefineSimpleKey{bib}{how}

\renewcommand{\eprint}[1]{#1}
\BibSpec{misc}{%
  +{}{\PrintAuthors}  {author}
  +{,}{ \textit}      {title}
  +{,}{ }             {how}
  +{}{ \parenthesize} {date}
  +{,} { available at \eprint}        {eprint}
  +{,}{ available at \url}{url}
  +{,}{ }             {note}
  +{.}{}              {transition}
}

\DeclareFontFamily{U}{matha}{\hyphenchar\font45}
\DeclareFontShape{U}{matha}{m}{n}{
      <5> <6> <7> <8> <9> <10> gen * matha
      <10.95> matha10 <12> <14.4> <17.28> <20.74> <24.88> matha12
      }{}
\DeclareSymbolFont{matha}{U}{matha}{m}{n}
\DeclareFontFamily{U}{mathb}{\hyphenchar\font45}
\DeclareFontShape{U}{mathb}{m}{n}{
      <5> <6> <7> <8> <9> <10> gen * mathb
      <10.95> mathb10 <12> <14.4> <17.28> <20.74> <24.88> mathb12
      }{}
\DeclareSymbolFont{mathb}{U}{mathb}{m}{n}
\DeclareMathSymbol{\ovoid}{3}{matha}{"6C}
\DeclareMathSymbol{\boxvoid}{2}{mathb}{"6C}

\mathchardef\mhyph="2D

\numberwithin{equation}{section}

\newtheorem{theorem}{Theorem}[section]
\newtheorem{corollary}[theorem]{Corollary}
\newtheorem{lemma}[theorem]{Lemma}
\newtheorem{proposition}[theorem]{Proposition}
\theoremstyle{remark}
\newtheorem{remark}[theorem]{Remark}

\theoremstyle{definition}
\newtheorem{definition}[theorem]{Definition}

\newcommand\bp{\begin{proof}}
\newcommand\ep{\end{proof}}
\newcommand\ee{\nopagebreak\mbox{\ }\hfill$\diamondsuit$}

\usepackage{wasysym}
\newcommand{\circt}%
{\mathbin{%
\mathchoice
{\ooalign{$\ocircle$\cr\hidewidth\raise-.15ex\hbox{$\scriptstyle\top\mkern2.05mu$}\cr}}
{\ooalign{$\ocircle$\cr\hidewidth\raise-.15ex\hbox{$\scriptstyle\top\mkern2.05mu$}\cr}}
{\ooalign{$\scriptstyle\ocircle$\cr\hidewidth\raise-.12ex\hbox{$\scriptscriptstyle\top\mkern1mu$}\cr}}
{\ooalign{$\scriptstyle\ocircle$\cr\hidewidth\raise-.12ex\hbox{$\scriptscriptstyle\top\mkern1mu$}\cr}}
}}

\DeclareMathOperator{\Ad}{Ad}

\DeclareMathOperator{\Irr}{Irr}
\DeclareMathOperator{\Rep}{Rep}

\DeclareMathOperator{\Tr}{Tr}

\newcommand{\C}{{\mathbb C}}
\newcommand{\FF}{{\mathbb F}}
\newcommand{\N}{{\mathbb N}}
\newcommand{\Z}{{\mathbb Z}}

\newcommand\T{{\mathbb T}}

\newcommand{\A}{{\mathcal A}}
\newcommand{\B}{{\mathcal B}}

\newcommand\CC{{\mathcal C}}

\newcommand{\F}{{\mathcal F}}
\newcommand\HH{\mathcal H}
\newcommand\I{\mathcal I}
\newcommand{\K}{{\mathcal K}}

\newcommand\OO{\mathcal O}

\newcommand\TT{\mathcal T}
\newcommand\U{\mathcal U}
\newcommand\eps{\varepsilon}

\begin{document}

\title[Subproduct systems]{Subproduct systems with quantum group symmetry}

\date{November 21, 2021; minor changes September 22, 2022}

\author{Erik Habbestad}
\address{University of Oslo, Mathematics institute}
\email{erikhab@math.uio.no}

\author{Sergey Neshveyev}
\email{sergeyn@math.uio.no}

\thanks{Supported by the NFR project 300837 ``Quantum Symmetry''.}

\begin{abstract}
We introduce a class of subproduct systems of finite dimensional Hilbert spaces whose fibers are defined by the Jones--Wenzl projections in Temperley--Lieb algebras. The quantum symmetries of a subclass of these systems are the free orthogonal quantum groups. For this subclass, we show that the corresponding Toeplitz algebras are nuclear C$^*$-algebras that are $KK$-equivalent to $\C$ and obtain a complete list of generators and relations for them. We also show that their gauge-invariant subalgebras coincide with the algebras of functions on the end compactifications of the duals of the free orthogonal quantum groups. Along the way we prove a few general results on equivariant subproduct systems, in particular, on the behavior of the Toeplitz and Cuntz--Pimsner algebras under monoidal equivalence of quantum symmetry groups.
\end{abstract}


\maketitle

\section*{Introduction}

The notion of a subproduct system, introduced by Shalit and Solel~\cite{MR2608451} and Bhat and Mukher\-jee~\cite{MR2646788}, lies at the intersection of two lines of research. One is dilation theory for semigroups of completely positive maps, the other is noncommutative function theory for row contractions. Recall that a row contraction is an $m$-tuple of operators $(S_1,\dots,S_m)$ such that $\sum^m_{i=1} S_iS_i^*\le1$. By a result of Popescu~\cite{MR1129595}, a universal, in a suitable sense, model for row contractions is provided by the creation operators $T_1,\dots T_m$ on the full Fock space~$\F(\C^m)$. If we are interested in the row contractions satisfying in addition algebraic relations $P_j(S_1,\dots,S_m)=0$ for some homogeneous polynomials $P_j$ in $m$ noncommuting variables, then such a model is often obtained by taking the compressions $S_i=e_\HH T_ie_\HH$ of the operators $T_i$ to the subspace $\F_\HH=\bigoplus^\infty_{n=0}H_n\subset\F(\C^m)$, where $H_n=I_n^\perp\subset(\C^m)^{\otimes n}$ and $I_n$ is the degree~$n$ component of the ideal $I\subset\C\langle X_1,\dots,X_m\rangle$ generated by $P_j$ \citelist{\cite{MR1668582}\cite{MR1749679}}. The collection $\HH=(H_n)^\infty_{n=0}$ has the property $H_{k+l}\subset H_k\otimes H_l$, which in the terminology of~\cite{MR2608451} means that it is a standard subproduct system of finite dimensional Hilbert spaces.

Our main object of interest in this paper is the structure of the Toeplitz algebra $\TT_\HH=C^*(S_1,\dots,S_m)\subset B(\F_\HH)$ and its quotient $\OO_\HH=\TT_\HH/\K(\F_\HH)$, called the Cuntz--Pimsner algebra of $\HH$ by Viselter~\cite{MR2949219}. Although the definition of these algebras is similar to such much more studied constructions as Cuntz--Krieger algebras, graph algebras and Pimsner algebras, a number of basic properties of the latter algebras do not have obvious analogues for $\TT_\HH$ and~$\OO_\HH$~\cite{MR2949219}, and a more or less complete understanding of $\TT_\HH$ and $\OO_\HH$ has been achieved only in a few cases. For example, Kakariadis and Shalit~\cite{MR3906397} made a comprehensive analysis of the case where the ideal $I$ is generated by monomials. In the present paper we consider principal ideals $I=\langle P\rangle$ generated by one quadratic polynomial. Furthermore, we require that the corresponding rank one projection $e=[\C P]\in B(\C^m\otimes\C^m)$ defines representations of Temperley--Lieb algebras $TL_n(\lambda^{-1})$ on $(\C^m)^{\otimes n}$ for some $\lambda\ge4$; recall that $TL_n(\lambda^{-1})$ is generated by projections $e_j$, $1\le j\le n-1$, satisfying the relations
$$
e_ie_j=e_je_i\quad\text{if}\quad |i-j|\ge2,\qquad e_ie_{i\pm1}e_i=\frac{1}{\lambda}e_i.
$$

Concretely, we consider the polynomials $P=\sum^m_{i,j=1}a_{ij}X_iX_j$ ($m\ge2$) such that $A\bar A$ is a scalar multiple of a unitary matrix, where $A=(a_{ij})_{i,j}$. We call such polynomials Temperley--Lieb. The simplest example of such a polynomial is $X_1X_2-X_2X_1$, which defines Arveson's symmetric subproduct system $SSP_2$~\cite{MR1668582}. In this case $\lambda=4$. Since the Temperley--Lieb algebras for $\lambda\ge4$ are all isomorphic, the entire collection that we consider can be thought of as a deformation/quantum analogue of $SSP_2$, and indeed, as we will see, the corresponding C$^*$-algebras $\TT_P$ and $\OO_P$ share a lot of properties. To be precise, our main results are proved for the subclass of polynomials such that $A\bar A=\pm1$. They have well-studied quantum symmetries~- the free orthogonal quantum groups~\cite{MR1382726}. The general case will be analyzed in detail in a separate publication.

Particular cases of the subproduct systems associated with Temperley--Lieb polynomials have already appeared in the literature. In addition to $SSP_2$, the case of polynomials $X_1X_2-qX_2X_1$ ($q\in\C^*$) has been studied in~\cite{MR1749679}, \cite{MR1905815} and~\cite{MR2608451}, although without a detailed analysis of the associated C$^*$-algebras. Our main class of examples, corresponding to the case $A\bar A=\pm1$, has been studied by Andersson~\cite{Anders}, with an emphasis on the gauge-invariant part of $\TT_P$. It should be said that his paper contains interesting ideas, but suffers from numerous imprecise and outright wrong constructions and arguments. Recently, the case of polynomials $\sum^m_{i=1}(-1)^iX_iX_{m-i+1}$ has been studied by Arici and Kaad~\cite{AK}, and some of our results strengthen and generalize the results in their paper.

\smallskip
The contents of the paper are as follows. In Section~\ref{sec:TL}, after a brief reminder on subproduct systems and associated C$^*$-algebras, we discuss the class of Temperley--Lieb polynomials. Using known formulas for the Jones--Wenzl projections in the Temperley--Lieb algebras we find some nontrivial relations in $\TT_P$.

In Section~\ref{sec:equivariant} we make a digression into the general theory of quantum group equivariant subproduct systems. By this we mean that each space $H_n$ is equipped with a unitary representation of a compact quantum group $G$ and the embedding maps $H_{k+l}\to H_k\otimes H_l$ are equivariant. A technical question that we study is under which conditions the action of $G$ on the Cuntz--Pimsner algebra $\OO_\HH$ is reduced, meaning that the averaging over $G$ defines a faithful cp map $\OO_\HH\to\OO^G_\HH$. Another question is what happens when we consider a monoidally equivalent compact quantum group $\tilde G$ and transform $\HH=(H_n)^\infty_{n=0}$ into a $\tilde G$-equivariant subproduct system~$\tilde\HH$. We show that the Toeplitz algebras $\TT_\HH$ and $\TT_{\tilde\HH}$ correspond to each other under the equivalence of categories of $G$- and $\tilde G$-C$^*$-algebras defined in~\cite{MR2664313}. The same is true for the reduced forms of the Cuntz--Pimsner algebras.

We want to make it clear that none of the results of this section is strictly speaking needed for the subsequent sections, but they provide a conceptual framework for the results in those sections.

In Section~\ref{sec:algebras} we return to the Temperley--Lieb polynomials and from this point onward concentrate on the polynomials $P=\sum^m_{i,j=1}a_{ij}X_iX_j$ such that $A\bar A=\pm1$. As we already mentioned, every such polynomial has a large quantum symmetry group $O^+_P$, where by large we mean that every $O^+_P$-module $H_n$ is irreducible. We show that by starting from Arveson's system $SSP_2$, corresponding to $P=X_1X_2-X_2X_1$, then moving to the polynomials $q^{-1/2}X_1X_2\pm q^{1/2}X_2X_1$ ($q>0$) and then considering the general case, the results of Section~\ref{sec:equivariant} quickly lead to the conclusion that the Cuntz--Pimsner algebra $\OO_P$ is $O^+_P$-equivariantly isomorphic to the linking algebra $B(SU_{\tau q}(2),O^+_P)$ (for suitable $q\in(0,1]$ and $\tau=\pm1$), which defines an equivalence between the representation categories of~$SU_{\tau q}(2)$ and $O^+_P$. We then find an explicit such isomorphism. It is worth stressing that the reason we know that $\OO_P$ cannot be smaller is that the action of~$O^+_P$ on $B(SU_{\tau q}(2),O^+_P)$ is reduced and ergodic, and therefore it does not admit any proper $O_P^+$-equivariant quotients. This is exactly the same argument as was used by Arveson for $SSP_2$ to conclude that the Cuntz--Pimsner algebra in his case is $C(S^3)$~\cite{MR1668582}. This property can be viewed as a replacement of the gauge-invariant uniqueness theorem, which fails for general subproduct systems~\cite{MR2949219}; see the recent paper by Dor-on~\cite{Dor-on} for a related discussion. Once $\OO_P$ has been computed, it is not difficult to show that the relations in $\TT_P$ that we found in Section~\ref{sec:TL} are complete.

In Section~\ref{sec:gauge} we consider the gauge-invariant part $\TT^{(0)}_P$ of $\TT_P$. Since $H_n$, $n\ge0$, exhaust the irreducible $O^+_P$-modules up to isomorphism, $\TT^{(0)}_P$ can be considered as an algebra of functions on a compactification of the dual discrete quantum group $\FF O_P$  of $O_P^+$. A compactification~$\overline{\FF O_P}$ of $\FF O_P$ has been constructed for the polynomials $P=q^{-1/2}X_1X_2\pm q^{1/2}X_2X_1$ by Tuset and the second author~\cite{MR2034922} using an analogy with the algebra of equivariant pseudo-differential operators of order zero on $SU(2)$, and for general $P$ by Vaes and Vergnioux~\cite{MR2355067} as a quantum analogue of the end compactification of a free group. We show that, excluding possibly the cases $O^+_P\cong SU_{\pm1}(2)$, we have $\TT^{(0)}_P=C(\overline{\FF O_P})$. As a consequence, we now have an explicit isomorphism $C(\partial{\FF O_P})\cong {}^\T B(SU_{\tau q}(2),O^+_P)$. The existence of such an isomorphism has been known, but in an indirect way, via an identification of both sides with the Martin boundary of~$\FF O_P$ \citelist{\cite{MR2034922}\cite{MR2355067}\cite{MR2400727}\cite{MR2664313}}.

In Section~\ref{sec:K-theory} we study $K$-theoretic properties of $\TT_P$. If we look again at $SSP_2$, then it follows already from an index theorem of Vegunopalkrishna~\cite{MR0315502} that the embedding map $\C\to \TT_{X_1X_2-X_2X_1}$ is a $KK$-equivalence. At this point it is natural to expect that the same is true for general $\TT_P$. For the polynomials $P=\sum^m_{i=1}(-1)^iX_iX_{m-i+1}$ this has indeed been shown by Arici and Kaad~\cite{AK}. They constructed an explicit inverse $\TT_P\to\C$ in the $KK$-category, which is a highly nontrivial task for $m\ge3$ (and for $m=2$ as well, if we want in addition $SU(2)$-equivariance). It is plausible that with some modification their construction and homotopy arguments work for general $P$. We take, however, a different route and use the Baum--Connes conjecture for~$\FF O_P$, as formulated by Meyer and Nest~\cite{MR2193334} and proved by Voigt~\cite{MR2803790}, to reduce the problem to comparing the $K$-theory of crossed products. We thus show that the embedding $\C\to\TT_P$ is a $KK^{O^+_P}$-equivalence without explicitly providing an inverse.

\smallskip\noindent
\textbf{Acknowledgements}: We are grateful to Stefaan Vaes and Christian Voigt for illuminating comments on linking algebras and $KK$-theory.

\bigskip

\section{Temperley--Lieb subproduct systems}\label{sec:TL}

Recall, following~\citelist{\cite{MR2608451}\cite{MR2646788}}, that a subproduct system  $\HH$ of finite dimensional Hilbert spaces (over the additive monoid $\Z_+$) is a sequence of Hilbert spaces $(H_n)^\infty_{n=0}$ together with isometries $w_{k,l}\colon H_{k+l}\to H_k\otimes H_l$ such that
$$
\dim H_0=1,\qquad \dim H_1=m<\infty,\qquad (w_{k,l}\otimes1)w_{k+l,n}=(1\otimes w_{l,n})w_{k,l+n}.
$$
By~\cite{MR2608451}*{Lemma~6.1}, we can assume that $H_0=\C$, $H_{k+l}\subset H_k\otimes H_l$ and the isometries $w_{k,l}$ are simply the embedding maps. The subproduct systems satisfying this stronger property are called standard. For such subproduct systems we denote by $f_n$ the projection $H_1^{\otimes n}\to H_n$.

\begin{remark}\label{rem:phase}
We could also start with formally weaker axioms for subproduct systems: it is enough to require that the identities $(w_{k,l}\otimes1)w_{k+l,n}=(1\otimes w_{l,n})w_{k,l+n}$ hold only up to phase factors. Indeed, the same arguments as in \cite{MR2608451}*{Lemma~6.1} show that we can still construct a standard subproduct system out of such a datum, and this standard system remains the same if we multiply~$w_{k,l}$ by phase factors. On a more conceptual level, the possibility of getting rid of phase factors in an essentially unique way follows from triviality of the cohomology groups $H^2(\Z_+;\T)$, $H^3(\Z_+;\T)$ (see~\cite{MR0077480}*{Proposition~X.4.1}). \ee
\end{remark}

Given a subproduct system $\HH=(H_n)^\infty_{n=0}$, the associated Fock space is defined by
$$
\F_\HH=\bigoplus^\infty_{n=0}H_n.
$$
For every $\xi\in H_1$, we define an operator
$$
S_\xi\colon\F_\HH\to\F_\HH\quad\text{by}\quad S_\xi\zeta=w_{1,n}^*(\xi\otimes\zeta)\quad\text{for}\quad \zeta\in H_n.
$$
We will usually fix an orthonormal basis $(\xi_i)^m_{i=1}$ in $H_1$ and write $S_i$ for $S_{\xi_i}$. The Toeplitz algebra~$\TT_\HH$ of $\HH$ is defined as the unital C$^*$-algebra generated by $S_1,\dots,S_m$.

If $\HH$ is standard and $H=H_1$, it is convenient to identify $\F_\HH$ with a subspace of the full Fock space
$$
\F(H)=\bigoplus^\infty_{n=0}H^{\otimes n}.
$$
Consider the operators $T_i\colon\F(H)\to\F(H)$, $T_i\zeta=\xi_i\otimes\zeta$, and the projection
$$
e_\HH\colon \F(H)\to\F_\HH.
$$
We obviously have $S_i=e_\HH T_i|_{\F_\HH}$.
Then also $S_i^*=e_\HH T_i^*|_{\F_\HH}$, but since $H_{n+1}\subset H_1\otimes H_n$, we actually have the stronger property
\begin{equation}\label{eq:S-vs-T}
S_i^*=T_i^*|_{\F_\HH}.
\end{equation}
Since $1-\sum^m_{i=1}T_iT_i^*$ is the projection onto $H^{\otimes 0}=\C$, we then get
\begin{equation}\label{eq:e-0}
e_0=1-\sum^m_{i=1}S_iS_i^*,
\end{equation}
where $e_0$ is the projection $\F_\HH\to H_0$. As the vacuum vector $\Omega=1\in H_0$ is cyclic for $\TT_\HH$, it follows that $\K(\F_\HH)\subset\TT_\HH$. The Cuntz--Pimsner algebra of $\HH$~\cite{MR2949219} is defined by
$$
\OO_\HH=\TT_\HH/\K(\F_\HH).
$$

Once we fix an orthonormal basis $(\xi_i)^m_{i=1}$ in $H$, it is convenient to identify the tensor algebra~$T(H)$ with the algebra $\C\langle X_1,\dots,X_m\rangle$ of polynomials in $m$ noncommuting variables. When we do this, we omit the symbol $\otimes$ for the product in $T(H)$, so we write $X_i\xi$ instead of $\xi_i\otimes\xi$.

By~\cite{MR2608451}*{Proposition~7.2}, there is a one-to-one correspondence between the standard subproduct systems with $H_1=H$ and the homogeneous ideals $I$ in $\C\langle X_1,\dots,X_m\rangle$ such that the degree one homogeneous component $I_1$ of $I$ is zero. Namely, given such an ideal $I$, we define
$$
H_n=I_n^\perp\subset H^{\otimes n}.
$$

In this work we mainly consider ideals $I=\langle P\rangle$ generated by one homogeneous polynomial of degree $2$. We denote by $\HH_P$ the corresponding standard subproduct system, but then normally use the subscript $P$ instead of $\HH_P$, so we write $\F_P$, $\TT_P$, etc.

\smallskip

We concentrate on polynomials having a particular symmetry.

\begin{definition}
Given a Hilbert space $H$ of finite dimension $m\ge2$, we say that a nonzero vector $\xi\in H\otimes H$ is \textbf{Temperley--Lieb}, if the corresponding projection $e=[\C\xi]\in B(H\otimes H)$ satisfies
\begin{equation}\label{eq:TL}
(e\otimes1)(1\otimes e)(e\otimes1)=\frac{1}{\lambda}e\otimes1\quad\text{in}\quad B(H\otimes H\otimes H)
\end{equation}
for some $\lambda>0$.
\end{definition}

As the following lemma shows, we then also have
$$
(1\otimes e)(e\otimes1)(1\otimes e)=\frac{1}{\lambda}1\otimes e.
$$

\begin{lemma}\label{lem:TL}
Assume $e_1$ and $e_2$ are projections in a C$^*$-algebra $A$ with a faithful tracial state~$\tau$ such that
$$
\tau(e_1)=\tau(e_2)\quad\text{and}\quad e_1e_2e_1=\frac{1}{\lambda}e_1\quad\text{for some}\quad \lambda>0.
$$
Then we also have $\displaystyle e_2e_1e_2=\frac{1}{\lambda}e_2$.
\end{lemma}

\bp
Consider the element $\displaystyle p=e_2-\lambda e_2e_1e_2$. Then $p^2=p$ and $\tau(p)=0$. Hence $p=0$.
\ep

In order to describe explicitly the Temperley--Lieb tensors, it is convenient to view the space $H\otimes H$ as the space of anti-linear operators $H\to H$, with $\zeta\otimes\eta$ corresponding to the operator $(\eta,\cdot)\zeta$. Equivalently, the tensor $\xi_A\in H\otimes H$ corresponding to an anti-linear operator $A\colon H\to H$ is given by
$$
\xi_A=\sum_i \xi_i\otimes A\xi_i,
$$
where $(\xi_i)_i$ is any orthonormal basis in $H$.

\begin{lemma}
Given $A\ne0$, the tensor $\xi_A\in H\otimes H$ is Temperley--Lieb if and only if $A^2$ is unitary up to a scalar factor, so that $(A^2)^*A^2=\alpha 1$ for some $\alpha>0$, and then $\lambda$ satisfying~\eqref{eq:TL} is given by $\lambda=\alpha^{-1}(\Tr A^*A)^2$.
\end{lemma}

\bp
Consider the polar decomposition $A=U|A|$. A priori $U$ is only an anti-linear partial isometry. We extend it to an anti-unitary $\tilde U$. Let $(\xi_i)_i$ be an orthonormal basis consisting of eigenvectors of $|A|$, $|A|\xi_i=\lambda_i \xi_i$. Put $\zeta_i=\tilde U \xi_i$. Let $(e_{ij})_{i,j}$ and $(f_{ij})_{i,j}$ be the matrix units in~$B(H)$ corresponding to $(\xi_i)_i$ and $(\zeta_i)_i$, resp. As $\xi_A=\sum_i \lambda_i\xi_i\otimes \zeta_i$, for the projection $e=[\C\xi_A]$ we have
$$
e=\beta\sum_{i,j}\lambda_i\lambda_j e_{ij}\otimes f_{ij},\quad\text{where}\quad\beta=(\Tr A^*A)^{-1}.
$$
Since $f_{i_1j_1}e_{i_2j_2}f_{i_3j_3}=(\xi_{i_2},\zeta_{j_1})(\zeta_{i_3},\xi_{j_2})f_{i_1j_3}$, we have
\begin{multline*}
(e\otimes1)(1\otimes e)(e\otimes1)\\
=\beta^3\sum_{i_1,i_2,i_3,j_1,j_2,j_3}\delta_{j_1,i_3}\lambda_{i_1}\lambda_{j_1}\lambda_{i_2}\lambda_{j_2}\lambda_{i_3}\lambda_{j_3}
(\xi_{i_2},\zeta_{j_1})(\zeta_{i_3},\xi_{j_2})e_{i_1j_3}\otimes f_{i_1j_3}\otimes f_{i_2j_2}.
\end{multline*}
It follows that~\eqref{eq:TL} holds if and only if the following identities are satisfied:
$$
\beta^3\lambda_{i_1}\lambda_{j_3}\sum_{i}\lambda_i^2\lambda_{i_2}\lambda_{j_2}(\xi_{i_2},\zeta_i)(\zeta_i,\xi_{j_2})=\delta_{i_2,j_2}\frac{\beta}{\lambda}\lambda_{i_1}\lambda_{j_3}.
$$
As $\lambda_{i_1}\lambda_{j_3}\ne0$ at least for some indices, it follows that if we introduce the unitary matrix $V=((\zeta_j,\xi_i))_{i,j}$ and the diagonal matrix $\Lambda=\operatorname{diag}(\lambda_1,\dots,\lambda_m)$, then we must have
$$
\lambda\beta^2\Lambda V\Lambda^2 V^*\Lambda=1.
$$
This is equivalent to invertibility of $\Lambda$ together with the identity $\lambda\beta^2 V\Lambda^2=\Lambda^{-2}V$. In terms of $|A|$ and $U$ this means that $\lambda\beta^2U|A|^2=|A|^{-2}U$, or equivalently, $\lambda\beta^2 AA^*=(A^*A)^{-1}$, which proves the lemma.
\ep

Using this lemma it is not difficult to parameterize the Temperley--Lieb tensors up to scalar factors and unitary transformations on $H$. Namely, since $(V\otimes V)\xi_A=\xi_{VAV^*}$ for any unitary~$V$ on~$H$, we want to classify up to unitary conjugacy antilinear operators $A\colon H\to H$ such that $A^2$ is unitary. If $A=U|A|$ is the polar decomposition of such an $A$, then $U|A|=|A|^{-1}U$. It follows that the spectrum $\sigma(|A|)$ of $|A|$ is closed under the transformation $\beta\mapsto\beta^{-1}$, and if $H_\beta\subset H$ is the spectral subspace for $|A|$ corresponding to $\beta$, then $UH_\beta=H_{\beta^{-1}}$. This implies that the collection of pairs $(\beta,Z_\beta)$, where $\beta\in S:=\sigma(|A|)\cap(0,1]$ and  $Z_\beta$ is the eigenvalues of $U^2$ on $H_\beta$ counting with multiplicity, is an invariant of the unitary conjugacy class of $A$. We claim that this is a complete invariant and the only restrictions on the sets $S$ and $Z_\beta$ are that $Z_1$ is closed under complex conjugation and
$$
\sum_{\beta\in S\cap(0,1)}2|Z_\beta|+|Z_1|=m.
$$

Indeed, it suffices to explain how to choose an orthonormal basis in which $A$ has a normal form. For every $\beta\in S\cap (0,1)$ we choose an orthonormal basis $(\xi_{\beta i})_i$ in $H_\beta$ consisting of eigenvectors of $U^2$, $U^2\xi_{\beta i}=z_{\beta i}\xi_{\beta i}$. Then the vectors $\zeta_{\beta i}=U\xi_{\beta i}$ form an orthonormal basis in~$H_{\beta^{-1}}$. Each two-dimensional space spanned by $\xi_{\beta i}$ and $\zeta_{\beta i}$ is invariant under $U$ and $|A|$, hence under $A$ and $A^*$, and the matrix of the restriction of $A$ to this space is
$$
\begin{pmatrix}
0 &  z_{\beta i}\beta^{-1}\\
\beta & 0
\end{pmatrix}.
$$

Next, assuming $1$ is an eigenvector of $|A|$, consider the space $H_1$. It is invariant under $U$, so here we can use the well-known description of anti-unitary matrices due to Wigner. It will be convenient to slightly modify it though. If $H_{1z}\subset H_1$ is the spectral subspace for $U^2$ corresponding to $z\in\T$, then $UH_{1z}=H_{1\bar z}$. It follows that, similarly to the previous paragraph, the direct sum of the spaces $H_{1z}$ for $z\ne\pm1$ decomposes into a direct sum of two-dimensional $U$-invariant subspaces such that on each of them $U$ has the form
$$
\begin{pmatrix}
0 &  z\\
1 & 0
\end{pmatrix},\quad 0<\arg(z)\le\pi.
$$
The same is true for $z=-1$, since if $\xi\in H_{1,-1}$, then $U\xi\perp\xi$, as
$$
(\xi,U\xi)=(UU\xi,U\xi)=-(\xi,U\xi).
$$

Finally, consider $H_{11}$. The real subspace of vectors fixed by $U$ is a real form of $H_{11}$. Choose an orthonormal basis $g_1,\dots,g_k$ in this Euclidean subspace. If $k=2l$ is even, we define an orthonormal basis $(\xi_j)_j$ in $H_{11}$ by
$$
\xi_j=\frac{1}{\sqrt{2}}(g_j+ig_{k-j+1})\quad\text{and}\quad \xi_{k-j+1}=\frac{1}{\sqrt{2}}(g_j-ig_{k-j+1}),\quad 1\le j\le l.
$$
If $k=2l+1$, then, in order to get a basis, in addition to the above vectors we need to take $\xi_{l+1}=g_{l+1}$. In both cases the matrix of $U|_{H_{11}}$ in the basis $(\xi_j)_j$ is
$$
\begin{pmatrix}
0 & & 1\\
 & \reflectbox{$\ddots$} & \\
 1 & & 0
\end{pmatrix}.
$$
This finishes the proof of our claim. The discussion can be summarized as follows.

\begin{proposition}\label{prop:TL-polynomials}
Assume $H$ is a Hilbert space of dimension $m=2l+r$ ($l\in\N$, $r\in\{0,1\}$) and $A\colon H\to H$ is an anti-linear operator such that $A^2$ is unitary. Then there is an orthonormal basis $(\xi_i)^m_{i=1}$ in which $A$ has the form
$$
\begin{pmatrix}
0 & & a_m\\
 & \reflectbox{$\ddots$} & \\
 a_1 & & 0
\end{pmatrix}\quad\text{for some}\quad a_i\in\C^*,\quad |a_ia_{m-i+1}|=1\quad\text{for all}\quad i,
$$
and then $\lambda$ satisfying~\eqref{eq:TL} for $\xi=\xi_A=\sum_i a_i\xi_i\otimes\xi_{m-i+1}$ is given by
$$
\lambda^{1/2}=\sum^m_{i=1}|a_i|^2\ge m.
$$
The numbers $a_1,\dots,a_m$ are uniquely determined by $A$ up to permutations of pairs of coefficients $(a_i,a_{m-i+1})$ ($1\le i\le l$) if we in addition require $0<a_i\le1$ ($1\le i\le l$), $a_{l+1}=1$ if $m$ is odd, and $0\le\arg(a_{m-i+1})\le\pi$ whenever $a_i=1$ ($1\le i\le l$).
\end{proposition}

Identifying $H$ with the space of homogeneous linear polynomials in $m$ variables, we can equivalently say that a homogeneous quadratic polynomial
$$
\sum^m_{i,j=1}a_{ij}X_iX_j
$$
is Temperley--Lieb if and only if the matrix $A=(a_{ij})_{i,j}$ is such that $A\bar A$ is a nonzero scalar multiple of a unitary matrix, where $\bar A=(\bar a_{ij})_{i,j}$. By a rescaling and a unitary change of variables, every such polynomial can be written as
$$
\sum^m_{i=1}a_iX_iX_{m-i+1},\quad\text{with}\quad |a_ia_{m-i+1}|=1\quad\text{for all}\quad i,
$$
and such a presentation is unique up to permutations of pairs $(a_i,a_{m-i+1})$ if we put restrictions on the coefficients $a_i$ as in Proposition~\ref{prop:TL-polynomials}.

\smallskip

Given a Temperley--Lieb polynomial $P=\sum^m_{i,j=1}a_{ij}X_iX_j$, consider the corresponding standard subproduct system $\HH_P=(H_n)^\infty_{n=0}$. If $e\in B(H\otimes H)$ is the projection onto $\C P$, then
$$
H_n=f_nH^{\otimes n},\qquad f_n=1-\bigvee^{n-2}_{i=0}1^{\otimes i}\otimes e\otimes1^{\otimes(n-i-2)}\quad\text{for}\quad n\ge2,
$$
and $f_0=1\in\C$, $f_1=1\in B(H)$. Since the projections $1^{\otimes i}\otimes e\otimes1^{\otimes(n-i-2)}$ satisfy the Temperley--Lieb relations, by computations of Jones~\cite{MR696688} and Wenzl~\cite{MR873400} (see \cite{neshveyev-tuset-book}*{Lemma~2.5.8} and its proof), we have
\begin{equation}\label{eq:Jones}
f_{n+1}=1\otimes f_{n}-[2]_q\phi(n)(1\otimes f_{n})(e\otimes 1^{\otimes(n-1)})(1\otimes f_{n}),\qquad n\ge0,
\end{equation}
where $q\in(0,1]$ is such that $\lambda^{1/2}=q+q^{-1}$, $\displaystyle [k]_q=\frac{q^k-q^{-k}}{q-q^{-1}}$ (with the convention $[k]_1=k$) and
\begin{equation}\label{eq:phi}
\phi(n)=\frac{[n]_{q}}{[n+1]_{q}}=\frac{q^{n}-q^{-n}}{q^{n+1}-q^{-n-1}}.
\end{equation}
Note that $\phi(n)\to q$ as $n\to+\infty$.

\begin{lemma}
For all $n\ge0$, we have $\dim H_n=[n+1]_t$, where $t\in(0,1]$ is such that $t+t^{-1}=m$.
\end{lemma}

\bp
It is known that the projection $[2]_q\phi(n)(1\otimes f_{n})(e\otimes 1^{\otimes(n-1)})(1\otimes f_{n})$ in~\eqref{eq:Jones} is equivalent to $e\otimes f_{n-1}$, see, e.g., again \cite{neshveyev-tuset-book}*{Lemma~2.5.8}. Hence the dimensions of the spaces $H_n=f_nH^{\otimes n}$ satisfy the recurrence relation
$$
\dim H_{n+1}=m \dim H_n-\dim H_{n-1}.
$$
This gives the result.
\ep

We next want to describe certain relations in $\TT_P$. As we will see later, these relations are complete, at least for Temperley--Lieb polynomials of a particular type.

\smallskip

Let us first introduce some notation. Since $\K(\F_P)\subset\TT_P$, we can view the algebra $c=C(\Z_+\cup\{+\infty\})$ of converging sequences as a subalgebra of $\TT_P$, with the projection $e_n\in c$ corresponding to the projection $\F_P\to H_n$. Denote by $\gamma$ the shift to the left on $c$.

\begin{proposition}\label{prop:relations}
Consider a polynomial $P=\sum^m_{i=1}a_iX_iX_{m-i+1}$ such that $|a_ia_{m-i+1}|=1$ for all $i$ ($m\ge2$). Let $q\in(0,1]$ be such that
$$
\sum^m_{i=1}|a_i|^2=q+q^{-1}.
$$
Then we have the following relations in $\TT_P$:
$$
fS_i =S_i\gamma(f)\quad (f\in c,\ 1\le i\le m),\qquad \sum^m_{i=1}S_iS_i^*=1-e_0,\qquad\sum^m_{i=1}a_iS_iS_{m-i+1}=0,
$$
$$
S_i^*S_j+a_i\bar a_j\phi S_{m-i+1}S_{m-j+1}^*=\delta_{ij}1\quad (1\le i,j\le m),
$$
where the element $\phi=(\phi(n))^\infty_{n=0}\in c$ is defined by~\eqref{eq:phi}.
\end{proposition}

\bp
The first relation simply reflects the fact that $S_iH_n\subset H_{n+1}$. The second relation is~\eqref{eq:e-0}, it holds in any subproduct system of finite dimensional Hilbert spaces. The third relation is an immediate consequence of the definition of $\HH_P$. So it is only the last relation that is not obvious.

Denote by $e_{ij}$ the matrix units in $B(H)$. The projection $e\in B(H\otimes H)$ onto $\C P$ is
$$
e=\frac{1}{[2]_q}\sum^m_{i,j=1}a_i\bar a_j e_{ij}\otimes e_{m-i+1,m-j+1}.
$$
Therefore by~\eqref{eq:Jones} we have
\begin{equation*}\label{eq:Jones2}
f_{n+1}=1\otimes f_{n}-\phi(n)\sum_{i,j}a_i\bar a_je_{ij}\otimes f_{n}(e_{m-i+1,m-j+1}\otimes 1^{\otimes(n-1)})f_{n}.
\end{equation*}
From this, for every $\xi\in H_n$, we get
$$
S_i^*S_j\xi=\delta_{ij}\xi-\phi(n)a_i\bar a_jf_{n}(e_{m-i+1,m-j+1}\otimes 1^{\otimes(n-1)})f_{n}\xi.
$$
By observing that by~\eqref{eq:S-vs-T} we have
$$
f_{n}(e_{kl}\otimes 1^{\otimes(n-1)})f_{n}=f_nT_kT_l^*f_n=S_kS_l^*
$$
on $H_n$, we obtain the last relation in $\TT_P$.
\ep

\bigskip

\section{Equivariant subproduct systems}\label{sec:equivariant}

As a preparation for a more thorough study of $\TT_P$, in this section we consider subproduct systems that are equivariant with respect to a compact quantum group.

\smallskip

Let us fix the notation and recall some basic notions, see~\cite{neshveyev-tuset-book} for more details. By a compact quantum group $G$ we mean a Hopf $*$-algebra $(\C[G],\Delta)$ that is generated by matrix coefficients of finite dimensional unitary right comodules. We have a one-to-one correspondence between such comodules and the finite dimensional unitary representations of $G$, that is, unitaries $U\in B(H_U)\otimes\C[G]$ such that $(\iota\otimes\Delta)(U)=U_{12}U_{13}$. Namely, the right comodule structure on $H_U$ is given by
$$
\delta_U\colon H_U\to H_U\otimes\C[G],\quad \delta_U(\xi)=U(\xi\otimes1).
$$
The tensor product $U\otimes V$ of two unitary representations (denoted also by $U\circt V$ or $U\times V$) is defined by $U_{13}V_{23}$. We denote by $\Irr(G)$ the set of the isomorphism classes of irreducible unitary representations of $G$. For every $s\in\Irr(G)$ we fix a representative $U_s\in B(H_s)\otimes\C[G]$.

Denote by $h$ the Haar state on $\C[G]$ and by $L^2(G)$ the corresponding GNS-space. We view $\C[G]$ as a subalgebra of $B(L^2(G))$. The (reduced) C$^*$-algebra $C(G)$ of continuous functions on~$G$ is defined as the norm closure of $\C[G]$.

Consider the $*$-algebra $\U(G)=\C[G]^*$ dual to the coalgebra $(\C[G],\Delta)$, with the involution $\omega^*(x)=\overline{\omega(S(x)^*)}$. Every finite dimensional unitary representation $U$ of $G$ defines a $*$-representation
$$
\pi_U\colon\U(G)\to B(H_U),\quad \pi_U(\omega)=(\iota\otimes\omega)(U).
$$
The representations $\pi_s$ defined by $U_s$, $s\in\Irr(G)$, give rise to an isomorphism
$$
\U(G)\cong \prod_{s\in\Irr(G)}B(H_s),
$$
and in what follows we are not going to distinguish between these two algebras. Consider the subalgebra
$$
c_c(\hat G)=\bigoplus_{s\in\Irr(G)}B(H_s)\subset\U(G).
$$
It coincides with the subalgebra of $C(G)^*\subset\U(G)$ spanned by the linear functionals $h(\cdot\,x)$ for $x\in\C[G]$. We remark that $c_c(\hat G)$ is weakly$^*$ dense in $C(G)^*$, since the state $h$ is faithful on~$C(G)$ and $\C[G]$ is dense in $C(G)$.

\smallskip

Assume now that $A$ is a C$^*$-algebra and $\alpha\colon A\to A\otimes C(G)$ is a $*$-homomorphism such that $(\alpha\otimes\iota)\alpha=(\iota\otimes\Delta)\alpha$.
We say that $\alpha$ is a (right) action of $G$ on $A$, or that $A$ is a $G$-C$^*$-algebra, if either of the following equivalent conditions is satisfied:
\begin{enumerate}
  \item the linear space $\alpha(A)(1\otimes C(G))$ is a dense subspace of $A\otimes C(G)$ (the Podle\'s condition);
  \item there is a dense $*$-subalgebra $\A\subset A$ on which $\alpha$ defines a coaction of the Hopf algebra $(\C[G],\Delta)$, that is, $\alpha(\A)\subset\A\otimes_{\mathrm{alg}}\C[G]$ and $(\iota\otimes\eps)\alpha(a)=a$ for all $a\in\A$,
\end{enumerate}
where $\otimes_{\mathrm{alg}}$ denotes the purely algebraic tensor product.
Then the largest subalgebra as in (ii) is given by
$$
\A=\operatorname{span}\{(\iota\otimes h)(\alpha(a)(1\otimes x)): a\in A,\ x\in\C[G]\},
$$
its elements are called regular. To put it differently, we have a left $c_c(\hat G)$-module structure on~$A$ defined by
$$
\omega\blacktriangleright a=(\iota\otimes\omega)\alpha(a),
$$
and then $\A=c_c(\hat G)\blacktriangleright A$.

An action is called reduced if $\alpha$ is injective, or equivalently, by faithfulness of $h$ on $C(G)$, if the conditional expectation
$$
E=(\iota\otimes h)\alpha\colon A\to A^G=\{a\in A:\alpha(a)=a\otimes1\}
$$
is faithful. For reduced actions the subalgebra of regular elements can also be described as
\begin{equation}\label{eq:regular}
\A=\{a\in A: \alpha(a)\in A\otimes_{\mathrm{alg}}\C[G]\},
\end{equation}
since if $\alpha(a)\in A\otimes_{\mathrm{alg}}\C[G]$, then $\alpha(a-(\iota\otimes\eps)\alpha(a))=0$.

Given any action $\alpha\colon A\to A\otimes C(G)$, we get a reduced action on $A_r=A/\ker\alpha$, since
$$
\ker\alpha=\{a\in A: E(a^*a)=0\}.
$$
We say that $A_r$ is the reduced form of the $G$-C$^*$-algebra $A$. Note that $A$ and $A_r$ have the same subalgebra $\A$ of regular elements.

If $G$ is coamenable, that is, the counit $\eps$ on $\C[G]\subset C(G)$ is bounded, then $(\iota\otimes\eps)\alpha=\iota$ on~$A$, so all actions of $G$ are reduced.

\begin{lemma}\label{lem:invariant-fd}
Assume $\alpha\colon A\to A\otimes C(G)$ is a reduced action and $X\subset A$ is a finite dimensional $G$-invariant subspace, so that $\alpha(X)\subset X\otimes C(G)$. Then all elements of $X$ are regular.
\end{lemma}

\bp
By assumption, $X$ is a $c_c(\hat G)$-submodule of $A$. As $X$ is finite dimensional, the action of $B(H_s)$ must be zero for all but a finite number of $s\in\Irr(G)$. In other words, there is a finite subset $F\subset\Irr(G)$ such that for all $\phi\in A^*$ the space $(\phi\otimes\iota)\alpha(X)\subset C(G)\subset L^2(G)$ is orthogonal to the matrix coefficients of $U_s$ for all $s\notin F$. It follows that this space is contained in the space $\C[G]_F$ spanned by the matrix coefficients of representations $U_s$, $s\in F$. But then $\alpha(X)\subset A\otimes \C[G]_F$, so $X$ consists of regular elements by~\eqref{eq:regular}.
\ep

Everything above of course also applies to the left actions $\alpha\colon A\to C(G)\otimes A$.

\smallskip

Assume now we are given two compact quantum groups $G$ and $\tilde G$. They are said to be monoidally equivalent if their categories of finite dimensional unitary representations are equivalent as C$^*$-tensor categories. By \citelist{\cite{MR1006625}\cite{MR2202309}}, every such equivalence, up to a natural isomorphism, is defined by a bi-Galois object/linking algebra. This is a unital C$^*$-algebra $B=B(G,\tilde G)$ equipped with two commuting reduced actions
$$
\delta\colon B\to C(G)\otimes B,\qquad \tilde\delta\colon B\to B\otimes C(\tilde G)
$$
that are ergodic and free; the precise meaning of the last property will not be important to us. Then an equivalence $F\colon \Rep G\to\Rep\tilde G$ is defined by mapping a finite dimensional unitary right comodule $(H_U,\delta_U)$ for $(\C[G],\Delta)$ into the cotensor product
$$
H_U\boxvoid_G B=\{X\in H_U\otimes B: (\delta_U\otimes\iota)(X)=(\iota\otimes\delta)(X)\},
$$
with the right $(\C[\tilde G],\Delta)$-comodule structure given by $\iota\otimes\tilde\delta$ and the scalar product
$$
(X,Y)=Y^*X\in {}^GB=\C1,
$$
where we interpret $\zeta^*\xi$ for $\xi,\zeta\in H_U$ as $(\xi,\zeta)$. The tensor structure on $F$ is defined by
$$
(H_U\boxvoid_G B)\otimes(H_V\boxvoid_G B)\to H_{U\otimes V}\boxvoid_G B,\quad X\otimes Y\mapsto X_{13}Y_{23}.
$$

In a similar way one can define an equivalence between the categories of reduced $G$-C$^*$-algebras and $\tilde G$-C$^*$-algebras, but there is a small technical issue that has to be taken care of~\cite{MR2664313}. Let $\B=\B(G,\tilde G)\subset B$ be the subalgebra of regular elements with respect to the $G$-action, equivalently, with respect to the $\tilde G$-action. Given a reduced action $\alpha\colon A\to A\otimes C(G)$, consider the subalgebra $\A\subset A$ of regular elements and define
$$
A\boxvoid_G B=\overline{\{X\in \A\otimes_{\mathrm{alg}}\B: (\alpha\otimes\iota)(X)=(\iota\otimes\delta)(X)\}}\subset A\otimes B.
$$
The right action of $\tilde G$ on $A\boxvoid_G B$ is again given by $\iota\otimes\tilde\delta$.

In some cases it is possible to use exactly the same definition of $A\boxvoid_G B$ as for finite dimensional comodules.

\begin{lemma}\label{lem:cotensor}
Assume $\alpha\colon A\to A\otimes C(G)$ is a reduced action of a compact quantum group $G$ on a C$^*$-algebra $A$  and there is a net of completely bounded $G$-equivariant maps $\theta_i\colon A\to A$ of finite rank such that $\theta_i(a)\to_i a$ for all $a\in A$ and $\sup_i\|\theta_i\|_{\mathrm{cb}}<\infty$. Then, for any compact quantum group $\tilde G$ monoidally equivalent to $G$, we have
$$
A\boxvoid_G B(G,\tilde G)=\{X\in A\otimes B(G,\tilde G): (\alpha\otimes\iota)(X)=(\iota\otimes\delta)(X)\}.
$$
\end{lemma}

\bp
If $(\alpha\otimes\iota)(X)=(\iota\otimes\delta)(X)$, then the elements $X_i=(\theta_i\otimes\iota)(X)$ have the same property and converge to $X$.
Therefore it suffices to show that $X_i\in\A\otimes_{\mathrm{alg}}\B(G,\tilde G)$. Lemma~\ref{lem:invariant-fd} implies that the finite dimensional space $\theta_i(A)$ is contained in $\A$, hence $X_i\in\A\otimes_{\mathrm{alg}}B(G,\tilde G)$ and then, using again that $(\alpha\otimes\iota)(X_i)=(\iota\otimes\delta)(X_i)$, we must have $X_i\in \A\otimes_{\mathrm{alg}}\B(G,\tilde G)$.
\ep

Consider a reduced action $\alpha\colon A\to A\otimes C(G)$ and a unitary, not necessarily finite dimensional, representation $U\in M(\K(H)\otimes C(G))$ of $G$. Assume we are given a covariant representation $\pi\colon A\to B(H)$, meaning that
$$
(\pi\otimes\iota)\alpha(a)=U(\pi(a)\otimes1)U^*.
$$
Given a monoidally equivalent compact quantum group $\tilde G$, we can extend the functor
$$
F=\cdot\,\boxvoid_G B(G,\tilde G)\colon\Rep G\to\Rep\tilde G
$$
to all unitary representations, since they decompose into irreducible ones. Therefore we get a Hilbert space $H\boxvoid_G B(G,\tilde G)$. By taking any state $\psi$ on $B(G,\tilde G)$ and considering the associated GNS-representation $\pi_\psi\colon B(G,\tilde G)\to B(H_\psi)$, we can view $H\boxvoid_G B(G,\tilde G)$ as a subspace of $H\otimes H_\psi$. It follows that $\pi\otimes\pi_\psi$ defines by restriction a representation $\pi\boxvoid\iota$ of $A\boxvoid_GB(G,\tilde G)$ on $H\boxvoid_G B(G,\tilde G)$.

\begin{lemma}\label{lem:injective-rep}
If the representation $\pi$ is faithful, then $\pi\boxvoid\iota$ is faithful as well.
\end{lemma}

\bp
Since the representation $\pi\boxvoid\iota$ is covariant and the $\tilde G$-action on $A\boxvoid_GB(G,\tilde G)$ is reduced, it suffices to check that the representation is faithful on
$$
(A\boxvoid_GB(G,\tilde G))^{\tilde G}=A^G\otimes1.
$$
Decompose $H$ into isotypical components, $H\cong\bigoplus_{s\in\Irr(G)}L_s\otimes H_s$,
where $L_s$ are some Hilbert spaces. Then the representation $\pi$ restricted to $A^G$ has the form
$\pi(a)=(\theta_s(a)\otimes1)_{s\in\Irr(G)}$ for some representations $\theta_s\colon A^G\to B(L_s)$. But then
$$
H\boxvoid_G B(G,\tilde G)\cong\bigoplus_{s\in\Irr(G)}L_s\otimes (H_s\boxvoid_G B(G,\tilde G)),
$$
and the representation $\pi\boxvoid\iota$ on $A^G\otimes1$ is simply
$$
(\pi\boxvoid\iota)(a\otimes 1)=(\theta_s(a)\otimes1\otimes1)_{s\in\Irr(G)},
$$
which is obviously faithful.
\ep

As we already mentioned, for the coamenable compact quantum groups all actions are reduced. As soon as $G$ is noncoamenable, there are nonreduced actions, for example, the action by right translations on the universal completion $C_u(G)$ of $\C[G]$. An interesting question is whether a quotient of a reduced action is reduced. We do not know an answer, but here is a partial result.

\begin{lemma}\label{lem:reduced}
Assume $\alpha\colon A\to A\otimes C(G)$ is a reduced action of a compact quantum group~$G$ on a C$^*$-algebra $A$, $I\subset A$ is a closed $G$-invariant ideal. Assume also that the following conditions are satisfied:
\begin{enumerate}
  \item $G$ is monoidally equivalent to a coamenable compact quantum group;
  \item there is a net of completely bounded $G$-equivariant maps $\theta_i\colon I\to I$ of finite rank such that $\theta_i(a)\to_i a$ for all $a\in I$ and $\sup_i\|\theta_i\|_{\mathrm{cb}}<\infty$.
\end{enumerate}
Then the action of $G$ on $A$ defines reduced actions of $G$ on $I$ and $A/I$.
\end{lemma}

For the proof we need the observation from~\cite{MR2202309} that if $\tilde G$ is coamenable, then $B(G,\tilde G)$ is a nuclear C$^*$-algebra. This follows from a general result in~\cite{MR1932666} stating that for coamenable~$\tilde G$ a $\tilde G$-C$^*$-algebra is nuclear if and only if the fixed point algebra is nuclear.

\bp[Proof of Lemma~\ref{lem:reduced}]
If $\A\subset A$ is the subalgebra of regular elements, then $\theta_i(I)\subset\A$ by Lemma~\ref{lem:invariant-fd}. As $\theta_i(a)\to_i a$ for all $a\in I$, it follows that $\I=I\cap\A$ is dense in $I$. Since $\alpha|_\I$  is a coaction of $(\C[G],\Delta)$ on $\I$, we conclude that $\alpha|_I$ is an action of $G$ on $I$. This action is obviously reduced.

Now, consider a coamenable compact quantum group $\tilde G$ monoidally equivalent to $G$.
Let us fix an equivalence of the representation categories of $G$ and $\tilde G$ and consider the corresponding bi-Galois objects $B(G,\tilde G)$ and $B(\tilde G,G)$. We then have bi-equivariant isomorphisms
$$
B(G,\tilde G)\boxvoid_{\tilde G}B(\tilde G,G)\cong C(G),\qquad B(\tilde G,G)\boxvoid_G B(G,\tilde G)\cong C(\tilde G),
$$
since the bi-Galois object corresponding to the identity functor on $\Rep G$ is $C(G)$. Then, modulo these isomorphisms, we have a natural isomorphism  $A\cong A\boxvoid_GB(G,\tilde G)\boxvoid_{\tilde G}B(\tilde G,G)$ for any reduced action $\alpha\colon A\to A\otimes C(G)$, defined by
$$
A\to A\boxvoid_G C(G),\quad a\mapsto\alpha(a).
$$

Put
$$
\tilde I=I\boxvoid_G B(G,\tilde G),\qquad \tilde A=A\boxvoid_G B(G,\tilde G), \qquad \tilde C=\tilde A/\tilde I.
$$
Denote by $\tilde\alpha$ the action of $\tilde G$ on $\tilde A$.
The action of $\tilde G$ on $\tilde C$ is reduced by coamenability. As $B(\tilde G,G)$ is a nuclear C$^*$-algebra, we have a short exact sequence
$$
0\to \tilde I\otimes B(\tilde G,G)\to \tilde A\otimes B(\tilde G,G)\to \tilde C\otimes B(\tilde G,G)\to0.
$$
By identifying $A$ with $\tilde A\boxvoid_{\tilde G} B(\tilde G,G)$ as explained above and denoting $\tilde C\boxvoid_{\tilde G} B(\tilde G,G)$ by $C$, we then get an exact sequence
$$
0\to (\tilde I\otimes B(\tilde G,G))\cap (\tilde A\boxvoid_{\tilde G} B(\tilde G,G))\to A\to C\to0.
$$

We claim that $(\tilde I\otimes B(\tilde G,G))\cap (\tilde A\boxvoid_{\tilde G} B(\tilde G,G))=I$. It is clear that
\begin{align*}
\tilde I\boxvoid_{\tilde G} B(\tilde G,G)&\subset (\tilde I\otimes B(\tilde G,G))\cap (\tilde A\boxvoid_{\tilde G} B(\tilde G,G))\\
&\subset\{x\in \tilde I\otimes B(\tilde G,G): (\tilde\alpha\otimes\iota)(x)=(\iota\otimes\tilde\delta)(x)\},
\end{align*}
where $\tilde\alpha$ and $\tilde\delta$ denote the actions of $\tilde G$ on $\tilde A$ and $B(\tilde G,G)$. The last set coincides with
$\tilde I\boxvoid_{\tilde G} B(\tilde G,G)=I$ by Lemma~\ref{lem:cotensor} applied to $\tilde G$ and the cb maps $\tilde\theta_i=\theta_i\otimes\iota|_{I\boxvoid_G B(G,\tilde G)}$ on $\tilde I$. This proves our claim.

Therefore we have a $G$-equivariant short exact sequence $0\to I\to A\to C\to0$ where all actions are reduced, which means that the action of $G$ on $A/I$ is reduced.
\ep

\begin{remark}\label{rem:exact}
What we tacitly used in the above argument is that if $A\to C$ is a surjective homomorphism of $G$-C$^*$-algebras, then we get a surjective map $\A\to\CC$ between the subalgebras of regular elements, and then,  for any compact quantum group $\tilde G$ monoidally equivalent to $G$, the map $\A\boxvoid_G\B(G,\tilde G)\to\CC\boxvoid_G\B(G,\tilde G)$ is surjective as well. This implies that if $0\to I\to A\to C\to0$ is a short exact sequence of reduced $G$-C$^*$-algebras, then the sequence of $\tilde G$-C$^*$-algebras
$$
0\to I\boxvoid_GB(G,\tilde G)\to A\boxvoid_GB(G,\tilde G)\to C\boxvoid_GB(G,\tilde G)\to0
$$
might not be exact in the middle, but at least the reduced form of $(A\boxvoid_GB(G,\tilde G))/(I\boxvoid_GB(G,\tilde G))$ is $C\boxvoid_GB(G,\tilde G)$.\ee
\end{remark}

Let us finally turn to equivariant subproduct systems. Assume $G$ is a compact quantum group and $\HH=(H_n)^\infty_{n=0}$ is a subproduct system of finite dimensional $G$-modules. By this we mean that we are given unitary representations $U_n$ of $G$ on $H_n$ and the structure maps $H_{k+l}\to H_k\otimes H_l$ are $G$-intertwiners.

In this case we have a unitary representation
$$
U_\HH=\bigoplus^\infty_{n=0} U_n\in M(\K(\F_\HH)\otimes C(G))
$$
of $G$ on $\F_\HH$. It defines reduced right actions of $G$ on $\K(\F_\HH)$ and $\TT_\HH$ by
$$
\alpha(T)=U_\HH(T\otimes1)U_\HH^*.
$$
To see that we indeed get an action on the Toeplitz algebra, fix an orthonormal basis in $H_1$ and let $u_{ij}$ be the matrix coefficients of $U_1$ in this basis. Then
\begin{equation}\label{eq:action-on-S}
\alpha(S_j)=\sum_i S_i\otimes u_{ij}.
\end{equation}
This implies that $\alpha$ defines a coaction of $(\C[G],\Delta)$ on the unital $*$-algebra generated by the operators $S_i$, hence it defines an action of $G$ on $\TT_\HH$.

We then get an action of $G$ on $\OO_\HH$. It is not clear to us whether this action is always reduced, but by Lemma~\ref{lem:reduced} we at least have the following.

\begin{proposition}\label{prop:reduced}
If $\HH=(H_n)^\infty_{n=0}$ is a subproduct system of finite dimensional $G$-modules and~$G$ is monoidally equivalent to a coamenable compact quantum group, then the action of~$G$ on~$\OO_\HH$ is reduced.
\end{proposition}

Assume now that $\tilde G$ is a compact quantum group monoidally equivalent to $G$ and fix a bi-Galois object $B=B(G,\tilde G)$. Given a subproduct system $\HH=(H_n)^\infty_{n=0}$ of finite dimensional $G$-modules, we get a subproduct system $\tilde\HH=(H_n\boxvoid _G B)^\infty_{n=0}$ of finite dimensional $\tilde G$-modules. Note that even if $\HH$ is standard, $\tilde\HH$ is usually not.

By definition we have a canonical unitary isomorphism $\F_{\tilde\HH}\cong\F_\HH\boxvoid_G B$. It follows that we get a representation of $\TT_\HH\boxvoid_G B$ on $\F_{\tilde\HH}$. This representation is faithful by Lemma~\ref{lem:injective-rep}.

\begin{proposition}\label{prop:monoidal-Toeplitz}
Assume $\HH=(H_n)^\infty_{n=0}$ is a subproduct system of finite dimensional $G$-modules and $\tilde G$ is monoidally equivalent to $G$. Consider the subproduct system $\tilde\HH=(H_n\boxvoid _G B(G,\tilde G))^\infty_{n=0}$. Then the canonical representation $\TT_\HH\boxvoid_G B(G,\tilde G)\to B(\F_{\tilde\HH})$ defines $\tilde G$-equivariant isomorphisms
$$
\K(\F_\HH)\boxvoid_G B(G,\tilde G)\cong\K(\F_{\tilde\HH}),\qquad \TT_\HH\boxvoid_G B(G,\tilde G)\cong\TT_{\tilde\HH}.
$$
\end{proposition}

\bp
The first isomorphism is a particular case of~\cite{MR2803790}*{Proposition~8.4}. The second will be obtained by similar arguments.

Since the representation $\TT_\HH\boxvoid_G B(G,\tilde G)\to B(\F_{\tilde\HH})$ is faithful, we can view $\TT_\HH\boxvoid_G B(G,\tilde G)$ as a subalgebra of $B(\F_{\tilde\HH})$. If $X=\sum_i\xi_i\otimes b_i\in \tilde H_1:=H_1\boxvoid_G B(G,\tilde G)$, then by definition the element $\sum_i S_{\xi_i}\otimes b_i\in \TT_\HH\boxvoid_G B(G,\tilde G)$ is the operator $S_X\in \TT_{\tilde\HH}$. Therefore $\TT_{\tilde\HH}\subset \TT_\HH\boxvoid_G B(G,\tilde G)$. It follows that
$$
\TT_{\tilde\HH}\boxvoid_{\tilde G} B(\tilde G,G)\subset \TT_\HH\boxvoid_G B(G,\tilde G)\boxvoid_{\tilde G} B(\tilde G,G)\subset B(\F_{\tilde\HH}\boxvoid_{\tilde G} B(\tilde G,G)).
$$
By identifying $\F_{\tilde\HH}\boxvoid_{\tilde G} B(\tilde G,G)$ with $\F_\HH$, this means that $\TT_{\tilde\HH}\boxvoid_{\tilde G} B(\tilde G,G)\subset\TT_\HH$. By swapping the roles of $G$ and $\tilde G$ we conclude that $\TT_{\tilde\HH}=\TT_\HH\boxvoid_G B(G,\tilde G)$ and $\TT_{\tilde\HH}\boxvoid_{\tilde G} B(\tilde G,G)=\TT_\HH$.
\ep

Together with Remark~\ref{rem:exact} this gives the following.

\begin{corollary}\label{cor:Cuntz-Pimsner}
If $\OO_{\HH,r}$ is the reduced form of the $G$-C$^*$-algebra $\OO_\HH$, then the $\tilde G$-C$^*$-algebra $\OO_{\HH,r}\boxvoid_G B(G,\tilde G)$ is isomorphic to the reduced form of $\OO_{\tilde\HH}$.
\end{corollary}

\bigskip

\section{Toeplitz and Cuntz--Pimsner algebras of Temperley--Lieb subproduct systems with large symmetry}\label{sec:algebras}

We now return to the subproduct systems defined by Temperley--Lieb polynomials and from now on consider only the polynomials $P=\sum^m_{i=1}a_iX_iX_{m-i+1}$ such that
$$
a_i\bar a_{m-i+1}=-\tau\in\{-1,1\}\quad\text{for all}\quad 1\le i\le m.
$$
A distinguishing property of these polynomials is that they have large quantum symmetries, the free orthogonal quantum groups~\cite{MR1382726}. Namely, for every $P$ as above consider the free orthogonal quantum group $O^+_P=O^+_{F_P}$ defined by the matrix
\begin{equation}\label{eq:F-matrix}
F_P=\begin{pmatrix}
0 & & a_m\\
 & \reflectbox{$\ddots$} & \\
 a_1 & & 0
\end{pmatrix}.
\end{equation}
The algebra $\C[O^+_P]$ of regular functions on $O^+_P$ is a universal unital $*$-algebra with generators~$u_{ij}$, $1\le i,j\le m$, and relations
$$
u_{ij}^*=a_i^{-1}a_ju_{m-i+1,m-j+1},\quad\text{the matrix}\quad(u_{ij})_{i,j}\quad\text{is unitary}.
$$
The comultiplication is given by $\Delta(u_{ij})=\sum_k u_{ik}\otimes u_{kj}$.

Consider the standard subproduct system $\HH_P=(H_n)^\infty_{n=0}$ defined by $P$. By definition, the quantum group $O^+_P$ has a unitary representation $U_1=(u_{ij})_{i,j}$ on $H_1$ and the representation $U_1\otimes U_1$ leaves $P\in H_1\otimes H_1$ invariant. Therefore we get unitary representations $U_n$ of~$O^+_P$ on $H_n=f_nH_1^{\otimes n}$. The spin~$\displaystyle\frac{n}{2}$ representations $U_n$, $n\ge0$, exhaust all irreducible representations of~$O^+_P$ up to equivalence, see~\cite{MR1378260} or \cite{neshveyev-tuset-book}*{Section~2.5}.

\begin{remark}
The quantum group $O_P^+$ is a genuine group only when $m=2$, $|a_1|=|a_2|=1$ and $\tau=1$, in which case it is $SU(2)$. The classical part of $O^+_P$, that is, the stabilizer of~$P$ in~$U(m)$, is generically rather small: if the numbers $|a_i|$ are all different, it consists of the unitary diagonal matrices $\operatorname{diag}(z_1,\dots z_m)$ such that $z_iz_{m-i+1}=1$.\ee
\end{remark}

Let us try to understand the Cuntz--Pimsner algebra $\OO_P$. In the simplest case $P=X_1X_2-X_2X_1$ the subproduct system $\HH_P$ is Arveson's symmetric subproduct system $SSP_2$ and $\OO_P$ is isomorphic to $C(S^3)\cong C(SU(2))$ by~\cite{MR1668582}*{Theorem~5.7}.

Next, take $q\in(0,1]$, $\tau=\pm1$ and consider the polynomial $P=q^{-1/2}X_1X_2-\tau q^{1/2}X_2X_1$. From Proposition~\ref{prop:relations} we see that the following relations are satisfied in $\OO_P$:
$$
S_1S_1^*+S_2S_2^*=1,\qquad q^{-1/2}S_1S_2-\tau q^{1/2}S_2S_1=0,
$$
$$
S_1^*S_1+S_2S_2^*=1,\qquad S_2^*S_2+q^2S_1S_1^*=1,\qquad S_1^*S_2-\tau qS_2S_1^*=0.
$$
These are the relations in $\C[SU_{\tau q}(2)]=\C[O^+_P]$ for the generators $u_{21}$ and $u_{22}$, which are usually denoted by $\gamma$ and $\alpha^*$, resp. In view of~\eqref{eq:action-on-S} and coamenability of $SU_{\tau q}(2)$ we conclude that there is a well-defined $SU_{\tau q}(2)$-equivariant surjective $*$-homomorphism $C(SU_{\tau q}(2))\to\OO_P$, $\gamma\mapsto S_1$, $\alpha\mapsto S_2^*$. This homomorphism must be injective, since the action of $SU_{\tau q}(2)$ on $C(SU_{\tau q}(2))$ defined by $\Delta$ is reduced and ergodic. Thus, $\OO_P\cong C(SU_{\tau q}(2))$, which is a quantum  analogue of Arveson's result for $q=\tau=1$.

Consider now the general case $P=\sum^m_{i=1}a_iX_iX_{m-i+1}$. Let $q\in(0,1]$ be such that
\begin{equation}\label{eq:q}
\sum^m_{i=1}|a_i|^2=q+q^{-1}.
\end{equation}
Then there is a unique up to a natural isomorphism monoidal equivalence between $O^+_P$ and $SU_{\tau q}(2)$, see \cite{MR2202309} or again \cite{neshveyev-tuset-book}*{Section~2.5}. Such an equivalence $F$ maps the spin~$\displaystyle\frac{n}{2}$ representations of $O^+_P$ into spin $\displaystyle\frac{n}{2}$ representations of $SU_{\tau q}(2)$. In both representation categories there is a unique up to a scalar factor morphism $H_{k+l}\to H_k\otimes H_l$. By Remark~\ref{rem:phase} it follows that $F$ maps~$\HH_P$ into an $SU_{\tau q}(2)$-equivariant subproduct system isomorphic to $\HH_{q^{-1/2}X_1X_2-\tau q^{1/2}X_2X_1}$. By Proposition~\ref{prop:reduced} and Corollary~\ref{cor:Cuntz-Pimsner} we conclude that
$$
\OO_P\cong C(SU_{\tau q}(2))\boxvoid_{SU_{\tau q}(2)}B(SU_{\tau q}(2),O^+_P)\cong B(SU_{\tau q}(2),O^+_P).
$$

Once we know that such an isomorphism is to be expected, it is not difficult to construct it from scratch using known generators and relations of $B(SU_{\tau q}(2),O^+_P)$. Namely, consider the matrices $F_P$, given by~\eqref{eq:F-matrix}, and
$$
F_{q,\tau}=\begin{pmatrix}
0 & -\tau q^{1/2}\\
q^{-1/2} & 0
\end{pmatrix}.
$$
Then by \cite{MR2202309}*{Theorem~5.5 and Remark~5.7}, $B(SU_{\tau q}(2),O^+_P)$ is a universal unital C$^*$-algebra with generators $y_{ij}$, $1\le i\le 2$, $1\le j\le m$, and relations
$$
Y=F_{q,\tau}\bar Y F_P^{-1},\qquad Y\quad\text{is unitary},
$$
where $Y=(y_{ij})_{i,j}$, $\bar Y=(y^*_{ij})_{i,j}$. The right action
$$
\delta\colon B(SU_{\tau q}(2),O^+_P)\to B(SU_{\tau q}(2),O^+_P)\otimes C(O^+_P)
$$
is given by $\displaystyle\delta(y_{ij})=\sum^m_{k=1}y_{ik}\otimes u_{kj}$, and the left action of $SU_{\tau q}(2)$ is defined in a similar way.

Consider the elements $y_j=y_{2j}$. The relation $Y=F_{q,\tau}\bar Y F_P^{-1}$ means that
$$
y_{1j}=q^{1/2}\bar a_jy_{m-j+1}^*,
$$
and then a simple computation shows that unitarity of $Y$ is equivalent to the relations
\begin{gather*}
\sum^m_{i=1}y_iy_i^*=1,\qquad\sum^m_{i=1}a_iy_iy_{m-i+1}=0,\\
y_i^*y_j+qa_i\bar a_jy_{m-i+1}y_{m-j+1}^*=\delta_{ij}1\quad (1\le i,j\le m).
\end{gather*}
These are exactly the relations in $\OO_P$ that we get from Proposition~\ref{prop:relations}. We thus arrive at the following result.

\begin{proposition}\label{prop:Cuntz-Pimsner}
Assume $P=\sum^m_{i=1}a_iX_iX_{m-i+1}$ is a polynomial such that $a_i\bar a_{m-i+1}=-\tau\in\{-1,1\}$ for all $i$ ($m\ge2$), and let $q\in(0,1]$ be given by~\eqref{eq:q}. Then we have an $O^+_P$-equivariant isomorphism $B(SU_{\tau q}(2),O^+_P)\cong\OO_P$, $y_i\mapsto S_i$ ($1\le i\le m$).
\end{proposition}

\bp
By the above discussion, the homomorphism $B(SU_{\tau q}(2),O^+_P)\to\OO_P$, $y_i\mapsto S_i$, is well-defined, surjective and $O^+_P$-equivariant. It must be injective, since the action of $O^+_P$ on the linking algebra $B(SU_{\tau q}(2),O^+_P)$ is reduced and ergodic.
\ep

Since $B(SU_{\tau q}(2),O^+_P)$ is nuclear by coamenability of $SU_{\tau q}(2)$, we get the following corollary.

\begin{corollary}
The C$^*$-algebras $\TT_P$ and $\OO_P$ are nuclear.
\end{corollary}

We are now ready to describe relations in $\TT_P$. Recall from Section~\ref{sec:TL} that we can view the C$^*$-algebra $c=C(\Z_+\cup\{+\infty\})$ of converging sequences as a subalgebra of $\TT_P$ and we denote by~$\gamma$ the shift to the left on $c$.

\begin{theorem}
Assume $P=\sum^m_{i=1}a_iX_iX_{m-i+1}$ is a polynomial such that $a_i\bar a_{m-i+1}=-\tau\in\{-1,1\}$ for all $i$ ($m\ge2$), and let $q\in(0,1]$ be given by~\eqref{eq:q}.  Then $\TT_P$ is a universal C$^*$-algebra generated by $S_1,\dots,S_m$ and $c$ satisfying the relations
$$
fS_i =S_i\gamma(f)\quad (f\in c,\ 1\le i\le m),\qquad \sum^m_{i=1}S_iS_i^*=1-e_0,\qquad\sum^m_{i=1}a_iS_iS_{m-i+1}=0,
$$
$$
S_i^*S_j+a_i\bar a_j\phi S_{m-i+1}S_{m-j+1}^*=\delta_{ij}1\quad (1\le i,j\le m),
$$
where the element $\phi\in c$ is defined by $\displaystyle \phi(n)=\frac{[n]_{q}}{[n+1]_{q}}$.
\end{theorem}

\bp
By Proposition~\ref{prop:relations} we already know that the above relations are satisfied in $\TT_P$. Consider a universal C$^*$-algebra $\tilde\TT_P$ generated $\tilde S,\dots,\tilde S_m$ and $c$ satisfying these relations. Let $\pi\colon\tilde\TT_P\to\TT_P$ be the quotient map.

Every element of the $*$-algebra generated by $\tilde S,\dots,\tilde S_m$ and $c$ can be written as a linear combination of elements of the form
$f \tilde S_{i_1}\dots \tilde S_{i_k}\tilde S^*_{j_1}\dots\tilde S^*_{j_l}$, where $f\in c$ and $k,l\ge0$. As $e_0\tilde S_i=0$, it follows $e_0\tilde\TT_Pe_0=\C e_0$, that is, $e_0$ is a minimal projection in $\tilde\TT_P$. Hence the closed ideal $I=\langle e_0\rangle$ generated by $e_0$ is isomorphic to the algebra of compact operators on some Hilbert space.

Next, since
$$
\sum^m_{i=1}\tilde S_ie_n\tilde S^*_i=e_{n+1}\sum^m_{i=1}\tilde S_i\tilde S^*_i=e_{n+1}(1-e_0)=e_{n+1},
$$
by induction on $n$ we conclude that $e_n\in I$ for all $n\ge0$, that is, $c_0\subset I$. It follows that the images of $\tilde S_i$ in $\tilde\TT_P/I$ satisfy the defining relations of $B(SU_{\tau q}(2),O^+_P)\cong\OO_P$. On the other hand, by construction, $\OO_P$ is a quotient of $\tilde\TT_P/I$. It follows that $\pi\colon\tilde\TT_P\to\TT_P$ gives rise to an isomorphism $\tilde\TT_P/I\cong\OO_P$.

Therefore we have a commutative diagram
$$
\xymatrix{0\ar[r] & I\ar[r]\ar[d]^{\pi|_I} & \tilde\TT_P\ar[r]\ar[d]^{\pi} & \OO_P\ar[r]\ar[d]^{\operatorname{id}} & 0\\
0\ar[r] & \K(\F_P)\ar[r] & \TT_P\ar[r] & \OO_P\ar[r] & 0}
$$
with exact rows and surjective $\pi$. As $I$ is a simple C$^*$-algebra, the map $\pi|_I$ must be an isomorphism, hence $\pi\colon\tilde\TT_P\to\TT_P$ is an isomorphism as well.
\ep

\begin{remark}
Since $\TT_P$ is generated by $S_1,\dots,S_m$, the relations in $\TT_P$ can be written without using the subalgebra~$c$. Namely, instead of the first two relations we could require the elements~$p_n$ defined by
$$
p_n=\sum^m_{i_1,\dots,i_n=1}S_{i_1}\dots S_{i_k}S^*_{i_k}\dots S^*_{i_1},\quad n\ge1,
$$
to be projections, and then define the element $\phi$ for the last relation by
$$
\phi=qp_1+\sum^\infty_{n=1}(\phi(n)-q)(p_n-p_{n+1}),
$$
which makes sense, as $p_1\ge p_2\ge\dots$ and therefore the projections $p_n-p_{n+1}$ are mutually orthogonal.

Indeed, then the projections $p_n$ together with the unit generate a copy of $c$, with $e_n=p_n-p_{n+1}$ for $n\ge0$, where $p_0=1$. As $p_{n+1}=\sum_i S_ip_nS_i^*$, we have
\begin{align*}
(1-p_{n+1})S_ip_nS_i^*(1-p_{n+1})&\le\sum_j(1-p_{n+1})S_jp_nS_j^*(1-p_{n+1})\\
&=(1-p_{n+1})p_{n+1}(1-p_{n+1})=0,
\end{align*}
so that $(1-p_{n+1})S_ip_n=0$. A similar computation gives $p_{n+1}S_i(1-p_n)=0$. It follows that $p_{n+1}S_i=S_ip_n$ for all $n\ge0$. This is equivalent to our first relation $f S_i=S_i\gamma(f)$.
\end{remark}

\bigskip

\section{Gauge action and compactifications of the dual discrete quantum groups}\label{sec:gauge}

We continue to consider the Temperley--Lieb polynomials $P=\sum^m_{i=1}a_iX_iX_{m-i+1}$ such that $a_i\bar a_{m-i+1}=-\tau\in\{-1,1\}$.

\smallskip

As the $O^+_P$-modules $H_n$ exhaust all irreducible representations of $O^+_P$ up to equivalence, the algebra of bounded functions on the dual discrete quantum group $\FF O_P=\widehat{O^+_P}$ is by definition
$$
\ell^\infty(\FF O_P)=\ell^\infty\text{-}\bigoplus^\infty_{n=0}B(H_n)\subset B(\F_P).
$$
For $x\in \ell^\infty(\FF O_P)$, we denote by $x_n$ its component in $B(H_n)$.

Consider the unitary representation $z\mapsto V_z$ of $\T$ on $\F_P$, where $V_z$ is the unitary that acts on~$H_n$ by the scalar $z^n$. Then we get an action $\Ad V$ of $\T$ on $B(\F_P)$, and $\ell^\infty(\FF O_P)$ coincides with the fixed point subalgebra of $B(\F_P)$ with respect to this action. The automorphisms $\Ad V_z$ leave~$\TT_P$ globally invariant and define the gauge action $\sigma$ of $\T$ on $\TT_P$, given by
$$
\sigma_z(S_i)=zS_i.
$$
Note that on $\OO_P\cong B(SU_{\tau q}(2),O^+_P)$ the gauge action coincides with the action of the maximal torus $\T\subset SU_{\tau q}(2)$. More precisely, the maximal torus is defined by the homomorphism
$$
\pi\colon C(SU_{\tau q}(2))\to C(\T),\quad \pi(u_{21})=0,\quad \pi(u_{22})(z)=\bar z.
$$
Therefore for the left action $\tilde\delta\colon B(SU_{\tau q}(2),O^+_P)\to C(SU_{\tau q}(2))\otimes B(SU_{\tau q}(2),O^+_P)$, $\tilde\delta(y_{ij})=\sum_k u_{ik}\otimes y_{kj}$, we have
$$
(\pi\otimes\iota)\tilde\delta(y_{2j})=(z\mapsto\bar z y_{2j})\in C(\T;B(SU_{\tau q}(2),O^+_P))=C(\T)\otimes B(SU_{\tau q}(2),O^+_P).
$$

Consider the gauge-invariant subalgebra
$$
\TT_P^{(0)}:=\TT_P^\sigma=\TT_P\cap \ell^\infty(\FF O_P)
$$
of $\TT_P$. Since it is unital and contains
$$
c_0(\FF O_P):=c_0\text{-}\bigoplus^\infty_{n=0}B(H_n)=\K(\F_P)^\sigma,
$$
it can be thought of as an algebra of continuous functions on a compactification of $\FF O_P$.

\smallskip

Assume $\sum_i|a_i|^2>2$. Following~\cite{MR2355067}, consider the inductive system of $O^+_P$-equivariant ucp maps
$$
\psi_{n,n+k}\colon B(H_n)\to B(H_{n+k}),\quad T\mapsto f_{n+k}(T\otimes1)f_{n+k},
$$
where we use that $H_{n+k}\subset H_n\otimes H_k$, and define
\begin{align*}
  C(\overline{\FF O_P}) & =\overline{\{x\in \ell^\infty(\FF O_P)\mid \psi_{n,n+k}(x_n)=x_{n+k}\ \text{for all}\ n\ \text{large enough and}\ k\ge0\}}\\
  & =\{x\in \ell^\infty(\FF O_P)\mid \lim_{n\to+\infty}\sup_{k\ge0}\|\psi_{n,n+k}(x_n)-x_{n+k}\|=0\}.
\end{align*}

Clearly, $c_0(\FF O_P)\subset C(\overline{\FF O_P})$. It is shown in~\cite{MR2355067} that $C(\overline{\FF O_P})$ is a C$^*$-algebra. This construction is a quantum analogue of the end compactification of a free group. Our goal in this section is to prove the following.

\begin{theorem}\label{thm:gauge-invariant-sub}
Assume $P=\sum^m_{i=1}a_iX_iX_{m-i+1}$ is a polynomial such that $a_i\bar a_{m-i+1}=-\tau\in\{-1,1\}$ for all $i$, and $\sum^m_{i=1}|a_i|^2>2$. Then $\TT_P^{(0)}=C(\overline{\FF O_P})$.
\end{theorem}

In the proof we are not going to use that $C(\overline{\FF O_P})$ is an algebra, so as a byproduct we will reprove that this is indeed the case.

\smallskip

Let $q\in(0,1)$ be given by~\eqref{eq:q}. The proof of the theorem relies on the following key known estimate.

\begin{lemma}\label{lem:VV}
There is a constant $C>0$ depending only on $q$ such that
$$
\|f_{n+1}-(1\otimes f_n)(f_n\otimes 1)\|\le Cq^n\quad\text{for all}\quad n\ge0.
$$
\end{lemma}

A stronger result is proved in~\cite{MR2355067}*{Lemma 8.4} using a generalization of Wenzl's recursive formula for $f_n$ established in~\cite{MR1446615}. Let us show that the particular case we need is a consequence of basic properties of the representation theory of $SU_q(2)$.

\bp[Proof of Lemma~\ref{lem:VV}]
It suffices to prove the estimate in the category $\Rep SU_{\tau q}(2)$ equivalent to $\Rep O_P^+$. Furthermore, since $\Rep SU_{-q}(2)$ can be obtained from $\Rep SU_q(2)$ by introducing new associativity morphisms such that $(H_1\otimes H_1)\otimes H_1\to H_1\otimes(H_1\otimes H_1)$ is the multiplication by~$-1$, it is enough to consider $\Rep SU_q(2)$. In other words, we may assume that $P=q^{-1/2}X_1X_2-q^{1/2}X_2X_1$.

As $H_{n+1}\subset (H_1\otimes H_n)\cap (H_n\otimes H_1)$, we just need to show that the restriction of $1\otimes f_n$ to $(H_n\otimes H_1)\ominus H_{n+1}$ has norm $\le Cq^n$. Since the $SU_q(2)$-module $(H_n\otimes H_1)\ominus H_{n+1}\cong H_{n-1}$ is irreducible, this norm is equal to $\|(1\otimes f_n)\xi\|$ for any unit vector $\xi$ in this module. As the vector~$\xi$ we take a highest weight vector. A formula for this vector is not difficult to find, up to a scalar factor it is
$$
[n]^{1/2}_qX_1^nX_2-q^{(n+1)/2}\zeta_nX_1,
$$
where $\zeta_n$ is a unit vector of weight $\displaystyle\frac{n}{2}-1$ in $H_n$, see the first paragraph of~\cite{neshveyev-tuset-book}*{Section~2.5}. Therefore we see that there is a unit vector $\xi_n\in (H_n\otimes H_1)\ominus H_{n+1}$ such that
$$
\|\xi_n-X_1^nX_2\|\le aq^n
$$
for a constant $a$ depending only on $q$. But then
$$
\|(1\otimes f_n)\xi_n\|\le aq^n+\|X_1 f_n(X_1^{n-1}X_2)\|\le aq^n+aq^{n-1}+\|f_n\xi_{n-1}\|=aq^n+aq^{n-1},
$$
which proves the lemma.
\ep

Similarly to the operators $S_i$ defined using the multiplication on the left, we can use the multiplication on the right and define
$$
R_i\colon\F_P\to\F_P,\quad R_i\xi=f_{n+1}(\xi X_i)\quad\text{for}\quad \xi\in H_n.
$$

\begin{lemma}\label{lem:commutators}
For all $1\le i,j\le m$, we have $[S_i,R_j]=0$. There is a constant $C$ depending only on $q$ such that
$$
\|[S_i^*,R_j]|_{H_n}\|\le Cq^n\quad\text{for all}\quad n\ge0.
$$
\end{lemma}

\bp
For $\xi\in H_n$, we have $S_iR_j\xi=f_{n+2}(X_i\xi X_j)=R_jS_i\xi$, which proves the first claim. To prove the second claim, take a unit vector $\xi\in H_n$ for some $n\ge1$. Then
$$
S_i^*R_j\xi=S_i^*f_{n+1}(\xi X_j)=T_i^*f_{n+1}(\xi X_j),
$$
where we  used~\eqref{eq:S-vs-T}. By Lemma~\ref{lem:VV} the last expression is $Cq^n$-close to
$$
T_i^*(1\otimes f_n)(\xi X_j)=f_n(T_i^*(\xi)X_j)=R_jS_i^*\xi,
$$
which gives the required estimate.
\ep

\begin{corollary}\label{cor:commutation-mod-compacts}
For every $S\in\TT_P$ and every $R$ in the C$^*$-algebra $C^*(R_1,\dots,R_m)$, we have
$[S,R]\in\K(\F_P)$.
\end{corollary}

\begin{remark}
The last corollary is true for $q=1$ as well, since in this case the proof of Lemma~\ref{lem:VV} gives
$\|f_{n+1}-(1\otimes f_n)(f_n\otimes 1)\|\le Cn^{-1/2}$.\ee
\end{remark}

Similarly to~\eqref{eq:e-0}, we have $\sum^m_{i=1}R_iR_i^*=1-e_0$. Define a contractive cp map
$$
\Theta\colon B(\F_P)\to B(\F_P),\quad \Theta(T)=\sum^m_{i=1}R_iTR_i^*.
$$
Since $\K(\F_P)\subset\TT_P$, by Corollary~\ref{cor:commutation-mod-compacts} this map leaves $\TT_P$ globally invariant. It also leaves $\ell^\infty(\FF O_P)$ globally invariant and we have
\begin{equation}\label{eq:theta}
\Theta^k(x)_{n+k}=\psi_{n,n+k}(x_n)\quad\text{for all}\quad x\in \ell^\infty(\FF O_P),\quad n,k\ge0.
\end{equation}

Denote by $\A_P$ the unital $*$-subalgebra of $\TT_P$ generated by the elements $S_i$, $1\le i\le m$. Let $\A^{(0)}_P$ be the gauge-invariant part of $\A_P$.

\begin{lemma}\label{lem:theta-estimate}
For every $x\in\A^{(0)}_P$, there exists a constant $C>0$ such that
$$
\|x_{n+k}-\Theta^k(x)_{n+k}\|\le Cq^n \quad\text{for all}\quad n,k\ge0.
$$
\end{lemma}

\bp
By Lemma~\ref{lem:commutators} we can find a constant $\tilde C$ depending on $x$ such that
$$
\|x_n-\Theta(x)_n\|\le\tilde Cq^n \quad\text{for all}\quad n\ge0.
$$
As $\Theta^l$ is a contraction mapping $B(H_{n+k-l})$ into $B(H_{n+k})$, it follows that
$$
\|\Theta^l(x)_{n+k}-\Theta^{l+1}(x)_{n+k}\|\le\tilde Cq^{n+k-l} \quad\text{for all}\quad n\ge0,\ k\ge l\ge0.
$$
Summing up over $l=0,\dots,k-1$ we get the required estimate, with $C=\tilde Cq(1-q)^{-1}$.
\ep

\bp[Proof of Theorem~\ref{thm:gauge-invariant-sub}]
The inclusion $\A^{(0)}_P\subset C(\overline{\FF O_P})$ follows from~\eqref{eq:theta} and Lemma~\ref{lem:theta-estimate}. Hence $\TT^{(0)}_P\subset C(\overline{\FF O_P})$.

In order to prove the opposite inclusion, take any $x\in\ell^\infty(\FF O_P)$ such that $$x_{n_0+k}=\psi_{n_0,n_0+k}(x_{n_0})$$ for some $n_0$ and all $k\ge0$. Since the maps $\psi_{n_0,n_0+k}$ are $O^+_P$-equivariant, we may assume that $x_{n_0}$ lies in a spin $l$ spectral component of $B(H_{n_0})$. Note that as $H_{n_0}\otimes H_{n_0}\cong H_0\oplus H_2\oplus\dots\oplus H_{2n_0}$, $l$ must be an integer $\le n_0$.

As $\OO_P^{(0)}\cong {}^\T B(SU_{\tau q}(2),O^+_P)$, the multiplicity of the spin $l$ component of $\OO_P^{(0)}$ is $1$. Fix a nonzero, hence injective, $O^+_P$-equivariant map $H_l\to \OO_P^{(0)}$ and lift it to an $O^+_P$-equivariant map $\pi\colon H_l\to\A^{(0)}_P$.
By Lemma~\ref{lem:theta-estimate} applied to $\pi(\xi)$ for the elements $\xi$ of a basis of $H_l$, we can find $C>0$ such that
\begin{equation}\label{eq:pi-theta}
\|\pi(\xi)_{n+k}-\Theta^k(\pi(\xi))_{n+k}\|\le Cq^n\|\xi\|\quad\text{for all}\quad \xi\in H_l\quad\text{and}\quad n,k\ge0.
\end{equation}

Next, we claim that there exist $\alpha>0$ and $n_1\ge0$ such that $\|\pi(\xi)_n\|\ge\alpha\|\xi\|$ for all $\xi\in H_l$ and $n\ge n_1$. If this is not true, then by compactness of the unit sphere of $H_l$ we can find a unit vector $\xi\in H_l$ such that $\liminf_n\|\pi(\xi)_n\|=0$. But then by~\eqref{eq:pi-theta} we must have $\lim_n\|\pi(\xi)_n\|=0$, that is, the image of~$\pi(\xi)$ in~$\OO^{(0)}_P$ is zero, which is a contradiction. Thus, our claim is proved.

Now, fix $\eps>0$ and choose $n\ge\max\{n_0,n_1\}$ such that $Cq^n\|x\|<\eps\alpha$. For the unique $\xi\in H_l$ such that $\pi(\xi)_n=x_n$ we have $\|\xi\|\le\alpha^{-1}\|x\|$. For all $k\ge0$ we have
$$
x_{n+k}=\Theta^k(x)_{n+k}=\Theta^k(\pi(\xi))_{n+k}.
$$
Hence, applying again~\eqref{eq:pi-theta}, we get
$$
\|x_{n+k}-\pi(\xi)_{n+k}\|\le Cq^n\|\xi\|<\eps.
$$
Therefore, modulo the compacts, $x$ is $\eps$-close to $\A^{(0)}_P$. Hence $x\in\TT^{(0)}_P$.
\ep

\bigskip

\section{K-theory}\label{sec:K-theory}

We continue to consider the same class of polynomials $P=\sum^m_{i=1}a_iX_iX_{m-i+1}$ as in the previous two sections.

\smallskip

Since $\OO_P\cong B(SU_{\tau q},O^+_P)$, by~\cite{MR3420332}*{Example~7.5} we have
$$
K_0(\OO_P)\cong\Z/(m-2)\Z,\qquad K_1(\OO_P)=\begin{cases}
                                              \Z, & \mbox{if}\quad m=2, \\
                                              0, & \mbox{if}\quad m\ge3.
                                            \end{cases}
$$
This already gives a lot of information about the $K$-theory of $\TT_P$ (for example, $K_1(\TT_P)=0$ if $m\ge3$), but is not quite enough to fully determine it. Our goal is to prove the following.

\begin{theorem}\label{thm:AK}
Assume $P=\sum^m_{i=1}a_iX_iX_{m-i+1}$ is a polynomial such that $a_i\bar a_{m-i+1}=-\tau\in\{-1,1\}$ for all $i$ ($m\ge2$). Then the embedding map $\C\to\TT_P$ is a $KK^{O^+_P}$-equivalence.
\end{theorem}

For the polynomials $P=\sum_i(-1)^iX_iX_{m-i+1}$ this upgrades the $KK^{SU(2)}$-equivalence of Arici--Kaad~\cite{AK} to the equivariant $KK$-category of the whole quantum symmetry group.

\smallskip

Our proof relies on the strong Baum--Connes conjecture for the dual of $SU_q(2)$ established by Voigt~\cite{MR2803790}. More precisely, we need the following standard consequence of this result.

\begin{proposition}
Assume $q\in[-1,1]\setminus\{0\}$, and $A$ and $B$ are separable $SU_q(2)$-C$^*$-algebras. Then a class $x\in KK^{SU_q(2)}(A,B)$ is a $KK^{SU_q(2)}$-equivalence if and only if it defines a $KK$-equivalence between $A\rtimes SU_q(2)$ and $B\rtimes SU_q(2)$.
\end{proposition}

\bp
A class $x$ is a $KK^{SU_q(2)}$-equivalence if and only if it induces an isomorphism
\begin{equation}\label{eq:KK-equiv}
KK^{SU_q(2)}(C,A)\cong KK^{SU_q(2)}(C,B)
\end{equation}
for every separable $SU_q(2)$-C$^*$-algebra $C$. The fact that $\widehat{SU_q(2)}$ is torsion-free and satisfies the strong Baum--Connes conjecture implies that the localizing subcategory of the triangulated category $KK^{SU_q(2)}$ generated by the separable C$^*$-algebras with trivial $SU_q(2)$-action coincides with $KK^{SU_q(2)}$~\cite{MR2803790}. Since the functors $KK^{SU_q(2)}(\cdot,A)$ and $KK^{SU_q(2)}(\cdot,B)$ are cohomological and transform the $c_0$-direct sums into direct products~\cite{MR2193334}, it follows that in order to check~\eqref{eq:KK-equiv} it suffices to consider~$C$ with trivial $SU_q(2)$-action. For such $C$ we have the Green--Julg isomorphisms
$$
KK^{SU_q(2)}(C,A)\cong KK(C,A\rtimes SU_q(2))\quad\text{and}\quad KK^{SU_q(2)}(C,B)\cong KK(C,B\rtimes SU_q(2)),
$$
see~\cite{vergnioux-thesis}*{Theorem~5.7}.
Therefore $x$ is a $KK^{SU_q(2)}$-equivalence if and only if it induces an isomorphism $KK(C,A\rtimes SU_q(2))\cong KK(C,B\rtimes SU_q(2))$ for every separable C$^*$-algebra $C$, that is, if and only if it defines a $KK$-equivalence between $A\rtimes SU_q(2)$ and $B\rtimes SU_q(2)$.
\ep

Combining this with the Universal Coefficient Theorem~\cite{MR894590}, we get the following.

\begin{corollary}\label{cor:UCT}
Assume $A$ and $B$ are separable $SU_q(2)$-C$^*$-algebras such that $A\rtimes SU_q(2)$ and $B\rtimes SU_q(2)$ satisfy the UCT. Then  a class $x\in KK^{SU_q(2)}(A,B)$ is a $KK^{SU_q(2)}$-equivalence if and only if it induces an isomorphism $K_*(A\rtimes SU_q(2))\cong K_*(B\rtimes SU_q(2))$.
\end{corollary}

Let us now fix our conventions and notation for the crossed products. Given a compact quantum group $G$, consider its right regular representation $V\in M(\K(L^2(G))\otimes C(G))$, defined~by
$$
V(a\xi_h\otimes\xi)=\Delta(a)(\xi_h\otimes\xi),\quad a\in C(G),\quad \xi\in L^2(G),
$$
where $\xi_h=1\in C(G)\subset L^2(G)$. Then $V(a\otimes1)V^*=\Delta(a)$. The integrated form of $V$ is the representation
$$
\rho=\pi_V\colon C^*(G)=c_0(\hat G)=c_0\text{-}\bigoplus_{s\in\Irr(G)}B(H_s)\to B(L^2(G)),
$$
$$
\rho(\omega)=(\iota\otimes\omega)(V)\quad\text{for}\quad\omega\in c_c(\hat G)\subset C^*(G).
$$
Explicitly, letting $\omega*a=(\iota\otimes\omega)\Delta(a)\in C(G)$ for $\omega\in c_c(\hat G)$ and $a\in C(G)$, we have
$$
\rho(\omega)a\xi_h=(\omega*a)\xi_h.
$$

Given an action $\alpha\colon A\to A\otimes C(G)$, the crossed product is defined by
$$
A\rtimes G=\overline{\alpha(A)(1\otimes\rho(C^*(G)))}\subset M(A\otimes\K(L^2(G))).
$$

\bp[Proof of Theorem~\ref{thm:AK}]
Let $q\in(0,1]$ be given by~\eqref{eq:q}.
By~\cite{MR2803790}*{Theorem~8.5}, the functor $A\mapsto A\boxvoid_{SU_{\tau q}(2)} B(SU_{\tau q}(2),O^+_P)$ extends to an equivalence of the equivariant $KK$-categories. Therefore by Proposition~\ref{prop:monoidal-Toeplitz} it suffices to prove the theorem for the polynomials
$$
P=q^{-1/2}X_1X_2-\tau q^{1/2}X_2X_1,\quad q\in(0,1],\quad \tau=\pm1.
$$

Consider the short exact sequence $0\to\K(\F_P)\to \TT_P\to \OO_P\to0$. Passing to crossed products we get a short exact sequence
\begin{equation}\label{eq:crossed1}
0\to\K(\F_P)\rtimes SU_{\tau q}(2)\to \TT_P\rtimes SU_{\tau q}(2)\to \OO_P\rtimes SU_{\tau q}(2)\to0.
\end{equation}
Since $\OO_P\cong C(SU_{\tau q}(2))$ by Proposition~\ref{prop:Cuntz-Pimsner}, by the Takesaki--Takai duality we have $ \OO_P\rtimes SU_{\tau q}(2)\cong\K(L^2(SU_{\tau q}(2)))$. Since the action of $SU_{\tau q}(2)$ on $\K(\F_P)$ is implemented by the unitary representation $U_P=\bigoplus^\infty_{n=0}U_n$, we also have an isomorphism
$$
\K(\F_P)\rtimes SU_{\tau q}(2)\cong \K(\F_P)\otimes\rho(C^*(SU_{\tau q}(2))),\quad X\mapsto U_P^*XU_P.
$$
As $C^*(SU_{\tau q}(2))=c_0\text{-}\bigoplus^\infty_{n=0}B(H_n)$, it follows that we can write~\eqref{eq:crossed1} as
\begin{equation}\label{eq:crossed2}
0\to c_0\text{-}\bigoplus^\infty_{n=0}\K(\F_P\otimes H_n)\to \TT_P\rtimes SU_{\tau q}(2)\to \K(L^2(SU_{\tau q}(2)))\to0.
\end{equation}
In this picture the canonical homomorphisms of $C^*(SU_{\tau q}(2))$ into $M(\K(\F_P)\rtimes SU_{\tau q}(2))$ and $\OO_P\rtimes SU_{\tau q}(2)$ are $(\pi_{U_P\otimes U_n})^\infty_{n=0}$
(since $(U_P^*)_{12}V_{23}(U_P)_{12}=(U_P)_{13}V_{23}$) and $\rho$, resp.

From~\eqref{eq:crossed2} it is clear that $\TT_P\rtimes SU_{\tau q}(2)$ is a type I C$^*$-algebra, hence it satisfies the UCT, and $K_1(\TT_P\rtimes SU_{\tau q}(2))=0$. Since $C^*(SU_{\tau q}(2))$ is also of type I with trivial $K_1$-group, by Corollary~\ref{cor:UCT} we just need to show that the canonical embedding $\pi\colon C^*(SU_{\tau q}(2))\to\TT_P\rtimes SU_{\tau q}(2)$ , $x\mapsto 1\otimes\rho(x)$, induces an isomorphism of the $K_0$-groups.

We identify $K_0(C^*(SU_{\tau q}(2)))$ with the representation ring $R(SU_{\tau q}(2))=\bigoplus^\infty_{n=0}\Z[U_n]$; in other words, $[U_n]\in K_0(C^*(SU_{\tau q}(2)))$ denotes the class of a rank one projection in $B(H_n)$. From~\eqref{eq:crossed2} we see that $K_0(\TT_P\rtimes SU_{\tau q}(2))$ is a free abelian group with generators $[p_n]$, $n\in\Z_+\cup\{\infty\}$,
where $p_n$ is a rank one projection in $\K(\F_P\otimes H_n)$ for $n\in\Z_+$ and $p_\infty$ is a projection in $\TT_P\rtimes SU_{\tau q}(2)$ such that its image in $\K(L^2(SU_{\tau q}(2)))$ is a rank one projection. As $p_\infty$ we can take the image of $1\in B(H_0)$ under $\pi\colon C^*(SU_{\tau q}(2))\to\TT_P\rtimes SU_{\tau q}(2)$, so that
$$
\pi_*([U_0])=[p_\infty].
$$

Now, let us fix $n\ge1$ and compute $\pi_*([U_n])$. We have
\begin{equation}\label{eq:K-map}
\pi_*([U_n])=c[p_\infty]+\sum^\infty_{k=0}c_k[p_k]
\end{equation}
for some $c,c_k\in\Z$, with only finitely many nonzero coefficients. The homomorphism $\TT_P\rtimes SU_{\tau q}(2)\to \K(L^2(SU_{\tau q}(2)))$ kills all the projections $p_n$, $n\ge0$. On the other hand, its composition with $\pi$ is the right regular representation $\rho$. As the multiplicity of $U_n$ in $V$ is $\dim H_n=n+1$, we conclude that $c=n+1$.

For $k\ge0$, consider the representation of $\TT_P\rtimes SU_{\tau q}(2)$ on $\F_P\otimes H_k$. Since the multiplicity of every isotopical component of $U_P\otimes U_k$ is finite, this is a representation by compact operators. Thus, we get a homomorphism $K_0(\TT_P\rtimes SU_{\tau q}(2))\to K_0(\K(\F_P\otimes H_k))$. Applying it to~\eqref{eq:K-map} we obtain
$$
(\pi_{U_P\otimes U_k})_*([U_n])=[n+1](\pi_{U_P\otimes U_k})_*([U_0])+c_k[p_k]\quad\text{in}\quad K_0(\K(\F_P\otimes H_k)).
$$
Therefore, if we denote by $m_{lk}$ the multiplicity of $U_l$ in $U_P\otimes U_k$, then
$$
c_k=m_{nk}-(n+1)m_{0k}.
$$
Since $U_P=\bigoplus^\infty_{i=0}U_i$ and $U_l\otimes U_k\cong U_{|l-k|}\oplus U_{|l-k|+2}\oplus\dots\oplus U_{l+k}$, we have
$$
m_{lk}=\frac{1}{2}(l+k-|l-k|)+1,\quad\text{hence}\quad c_k=\frac{1}{2}(k-n-|k-n|).
$$

To summarize, for all $n\ge0$ we have
$$
\pi_*([U_n])=(n+1)[p_\infty]+\sum^{n-1}_{k=0}(k-n)[p_k].
$$
This shows that $\pi_*\colon R(SU_{\tau q}(2))\to K_0(\TT_P\rtimes SU_{\tau q}(2))$ is indeed an isomorphism.
\ep

\bigskip\bigskip

\bigskip

\end{document}